\documentclass[11pt,reqno]{amsart}
\usepackage{amsbsy,amsfonts,amsmath,amssymb,amscd,amsthm}
\usepackage{mathrsfs}
\usepackage{subfigure,enumitem,color,graphicx}
\usepackage[a4paper,text={6.1in,8in},centering,includefoot,foot=1in]{geometry}
\usepackage{pxfonts}
\usepackage[colorlinks=true,linkcolor=blue,citecolor=green]{hyperref}
\usepackage{srcltx}

\small\normalsize
\theoremstyle{plain}
\newtheorem{maintheorem}{Theorem}
\newtheorem{thm}{Theorem}[section]
\newtheorem{cor}[thm]{Corollary}
\newtheorem{lem}[thm]{Lemma}
\newtheorem{claim}[thm]{Claim}
\newtheorem{prop}[thm]{Proposition}
\newtheorem{defi}[thm]{Definition}
\theoremstyle{definition}
\newtheorem{rem}[thm]{Remark}

\numberwithin{equation}{section}

\newcommand{\eqdef}{\stackrel{\scriptscriptstyle\rm def}{=}}
\DeclareMathOperator{\IFS}{IFS}


\begin{document}

\title{Symbolic blender-horseshoes and applications}
\author[Barrientos, Ki \& Raibekas]{Pablo G. Barrientos, Yuri Ki and Artem Raibekas \vspace{1cm}}

\address{
\footnotesize
Instituto de Matem\'atica e Estat\'istica,
Universidade Federal Fluminense,\newline
Rua M\'ario Santos Braga s/n - Campus Valonguinhos, Niter\'oi, RJ, Brasil
}
\email{barrientos@uniovi.es}

\address{
\footnotesize
Instituto de Ci\^encias Exatas, Universidade Federal Fluninense, \newline
Rua Desembargador Ellis Hermydio Figueira, 783, Bairro Aterrado, Volta Redonda, RJ, Brasil
}
\email{yuriki@impa.br}

\address{
\footnotesize
Instituto de Matem\'atica e Estat\'istica,
Universidade Federal Fluminense,\newline
Rua M\'ario Santos Braga s/n - Campus Valonguinhos, Niter\'oi, RJ, Brasil
}
\email{artem@impa.br}

\begin{abstract}
We study partially-hyperbolic skew-product maps over the Bernoulli shift with
H\"older dependence on the base points. In the case of contracting fiber maps,
symbolic blender-horseshoe is defined as an invariant set
which meets any almost horizontal disk in a robust sense. These invariant sets are understood as blenders with center stable bundle of any dimension. We then give necessary conditions (covering property) on an iterated function system such that the relevant skew-product has a symbolic blender-horseshoe. We use this local plug to yield robustly non-hyperbolic transitive diffeomorphisms and robust heterodimensional cycles of co-index equal to the dimension of the central
direction.
\end{abstract}
~\vspace{-1cm} \maketitle


\section{Introduction}
\label{s.intro}

Currently hyperbolic systems are quite well understood from the
topological and statistical perspectives, see~\cite{BDV05}. A
natural objective is to characterize obstructions to
hyperbolicity. In~\cite{P00}, Palis conjectured that generation of
cycles (heterodimensional cycles and homoclinic tangencies) is the
only obstruction: every diffeomorphism can be approximated
either by a hyperbolic one or by one with a cycle. This conjecture
was proved for surface $C^1$-diffeomorphisms in~\cite{PS00} and in higher dimensions  there are some
partial results, see \cite{CP10}.

On the other hand,  dynamical properties
which exhibit some ``persistence''  play an important role in dynamics. A property of a diffeomorphism is robust  if there is a neighborhood consisting of  diffeomorphisms
satisfying such a property. Bearing in mind this concept,
Bonatti formulated in~\cite{Bonatti2004} a strong version of ``hyperbolicity versus cycles'' conjecture: the set of hyperbolic diffeomorphisms and the set of
diffeomorphisms with robust cycles are two  open sets whose union
is dense in the space of $C^1$-diffeomorphisms. This conjecture
was proved in~\cite{BD08} for the so-called ``tame systems''.

In this paper we study the generation of $C^1$-robust cycles continuing the work of~\cite{BD08,BD11}. One of the main tools used are the so-called ``blenders'', which are hyperbolic sets in some sense similar to the thick
horseshoes introduced by Newhouse~\cite{New70} in the setting of homoclinic tangencies.
Roughly, a blender is a local plug that guarantees  the  closure of an invariant  manifold of
a hyperbolic set to be greater than the dimension of its unstable bundle.

Besides playing a key role in the generation of robust cycles (see~\cite{BD08,BD11}), blenders are also important in other settings such as construction of non-hyperbolic robustly transitive sets~\cite{BD96}, discontinuity of the dimension
of hyperbolic sets~\cite{BDV95},  stable ergodicity~\cite{HHTU}, and Arnold difusion in symplectic and Hamiltonian settings~\cite{NP11, FonMar}.

The papers~\cite{BD96, BD08, BD11} consider blenders having  one-dimensional ``central'' direction, see Definition~\ref{d.blender}.
Here, following the ideas in \cite{NP11}  we study blenders with ``central'' direction
of any dimension, which enables us to consider these sets in  more general settings, and
as applications  we get the next results.

\subsection{Robust heterodimensional cycles and mixing sets.}
\label{s.aplica}
Let $M$ be a manifold and consider $\mathrm{Diff}^1(M)$ the set of diffeomorphisms $f : M \to M$ endowed with the $C^1$-topology.
A \emph{heterodimensional cycle} for 
$f \in \mathrm{Diff}^1(M)$
between a pair of hyperbolic transitive sets
$\Omega_1$ and $\Omega_2$
of different \emph{indices}
(dimension of their stable bundles) is defined by the relations
$$
W^s(\Omega_1,f)\cap W^u(\Omega_2,f)\ne\emptyset
\quad \mbox{and} \quad
W^u(\Omega_1,f)\cap W^s(\Omega_2,f)\ne\emptyset.
$$
The {\emph{co-index}} 
is
the absolute value of the difference between the
indices of $\Omega_1$,
$\Omega_2$.
We are interested in cycles that cannot be destroyed by
perturbations.

\begin{defi}[Robust heterodimensional cycle]
\label{d.ciclo robusto}
A diffeomorphism $f$ has a \emph{$C^1$-robust heterodimensional cycle}
associated with its transitive hyperbolic sets $\Omega_1$ and $\Omega_2$
if there is a $C^1$-neighborhood $\mathcal{U}$ of
$f$ such that every $g \in \mathcal{U}$ has a heterodimensional cycle
associated with the hyperbolic continuations $\Omega_1^g$ and $\Omega_2^g$
of $\Omega_1$ and $\Omega_2$, respectively.
\end{defi}

In \cite{BD08} it is proved that  cycles of co-index
one 
yield $C^1$-robust cycles after small perturbations.
A key ingredient is the construction of blenders
with one dimensional center-stable direction.
Here we consider $F \colon N \to N$ a $C^1$-diffeomorphism with a horseshoe,
$M$ a manifold of dimension $c\geq 1$ and $\mathrm{id}$ the identity map.
We prove that there are maps arbitrarily close to $F\times \mathrm{id}$
having robust  cycles.
The following result  is a generalization of
\cite[Theorem 2.4]{BD08}.

\begin{maintheorem}
\label{t.ciclo robusto}
Let $F:N\to N$ be a $C^1$-diffeomorphism with a Smale horseshoe.
Then there is an arc
$\{ f_\varepsilon \}_{\varepsilon \in (0, \varepsilon_0]}$ of $C^1$-diffeomorphisms
on $N\times M$ such that $f_0=F\times \mathrm{id}$ and $f_\varepsilon$ has
a $C^1$-robust heterodimensional cycle of co-index~$c$
for all $0<\varepsilon<\varepsilon_0$.
\end{maintheorem}




The construction in \cite{BD96} of persistent non-hyperbolic transitive diffeomorphisms used a ``chain of blenders'',
which is a connected family of blenders.  
Here, we provide an alternative construction using {\it symbolic blenders} (the precise definition will be given in
Section~\ref{ss.symbolicblenderhorseshoes}),
thus obtaining a slightly stronger, topologically mixing, version of the theorem.
An invariant set $X$ is {\it topologically mixing} for a map $g$ if for any pair of open sets $U$, $V$ of $X$ there is a $n_0>0$ such that $g^{n}(U) \cap V \neq \emptyset$ for all $n\geq n_0$.

\begin{maintheorem}
\label{t.app mix} Let $F:N\to N$ be a $C^1$-diffeomorphism with a
Smale horseshoe $\Lambda$. Then there is an arc
$\{f_\varepsilon\}_{\varepsilon \in (0,\varepsilon_0]}$ of
$C^1$-diffeomorphisms on $N\times M$ such that $f_0=F\times
\mathrm{id}$ and for every small enough $C^1$-perturbation $g$ of $f_\varepsilon$, $0<\varepsilon<\varepsilon_0$, there is a $g$-invariant set $\Delta_g \subset N\times M$ homeomorphic to $\Lambda\times M$ which is a  non-hyperbolic
topologically
mixing set for $g$.
\end{maintheorem}

\subsection{The set of symbolic skew-products}
Blenders appear in a natural way in the context  of skew-product diffeomorphisms
$$
 f: [-1,1]^{n}\times [-1,1]
\to \mathbb{R}^n \times \mathbb{R},\qquad f(x,y)=
\big( F(x),\phi(x,y) \big), \vspace{0.2cm}
$$
where $n\geq 2$, $F\colon\mathbb{R}^n\to \mathbb{R}^n$ has a horseshoe
$\Lambda \subset [-1,1]^n$, and the fiber maps $\phi(x,\cdot)\colon
\mathbb{R} \to \mathbb{R}$ are contractions.
Since $F|_\Lambda$ is conjugate to a shift, this map can be studied symbolically. This leads to consider {\emph{symbolic skew-product maps}} $\Phi$ over the
shift map $\tau:\Sigma_k\eqdef \{1,\ldots,k\}^\mathbb{Z}\to \Sigma_k$ of $k\geq 2$ symbols of the
form
\begin{equation}
 \label{e:sym-skew}
   \Phi: \Sigma_k \times M  \to \Sigma_k \times M, \qquad
   \Phi(\xi,x)= \big( \tau(\xi), \phi_\xi(x) \big),
\vspace{0.2cm}
\end{equation}
where $M$ is a manifold and $\phi_\xi: M
\to M$ are diffeomorphisms depending continuously with respect to $\xi$. The set $\Sigma_k$ is
called \emph{base} and $M$ is the \emph{fiber}.
In this paper we will work with maps of the form~\eqref{e:sym-skew} and in order to emphasize the role of the fiber maps we write $\Phi=\tau\ltimes\phi_\xi$.

\label{ss:thesetS}
Let $D$ be a bounded open set of $M$.
Given positive constants $0<\lambda<\beta$, a map $\phi\colon
\overline D\to D$ is said to be a \emph{$(\lambda,
\beta)$-Lipschitz} on $\overline D$ if
$$
\lambda\,\|x-y\|< \|\phi (x)- \phi (y)\| < \beta \, \|x-y\|, \quad
\mbox{ for all $x, y \in \overline D$,}
$$
where $\|x-y\|$ denotes the distance between $x$ and $y$ in $M$.
Fix $0<\nu<1$ and the metric
\begin{equation*}
\label{e:metrica}
d_{\Sigma_k}(\xi,\zeta)\eqdef \nu^{\ell}, \ \ \ell=\min\{i\in
\mathbb{Z}^+: \xi_i\not=\zeta_i\ \text{or} \
\xi_{-i}\not=\zeta_{-i} \}.
\end{equation*}
For a fixed $0\leq \alpha \leq 1$, we say that $\Phi=\tau\ltimes\phi_\xi$ is
\emph{locally $\alpha$-H\"older continuous} (or that its fiber maps $\phi_\xi$ depend
\emph{locally $\alpha$-H\"older continuously} on $\overline{D}$ with respect to
the base point $\xi$) if there is a constant $C\ge 0$ such
that
\begin{equation}
\label{eq:Holder}
\|\phi^{\pm 1}_\xi(x)-\phi^{\pm 1}_\zeta(x)\|
\leq C \, d_{\Sigma_k}(\xi,\zeta)^{\alpha}, \quad
\mbox{for all $x\in \overline D$} 
\end{equation}
and every $\xi, \zeta \in \Sigma_k$ with $\xi_0=\zeta_0$.
We denote by $C_\Phi$ the smallest non-negative constant
satisfying~\eqref{eq:Holder} and call it the \emph{local H\"older
constant} of $\Phi$ on $\overline{D}$.

\begin{defi}[Sets of symbolic skew-products]
\label{d.thesetS}
Let $D\subset M$ be a bounded open set,
$r\geq 0$, $0<\lambda<\beta$, $0\leq \alpha \leq 1$ and $k>1$.
We define $\mathcal{S}^{r,\,\alpha}_{k,\lambda,\beta}(D)$
as the \emph{set of symbolic skew-product maps} $\Phi =\tau
\ltimes \phi_\xi$ as in~\eqref{e:sym-skew}
such that
\begin{itemize}
\item $\phi_\xi$ is $C^r$-$(\lambda,\beta)$-Lipschitz on $\overline{D}$ for all $\xi\in\Sigma_k$,
\item $\phi_\xi$ depends locally $\alpha$-H\"older continuously on $\overline{D}$  with respect to $\xi$.
\end{itemize}
The set
 $\mathcal{S}^{r,\alpha}_{k,\lambda,\beta}(D)$ is endowed
with the distance
\begin{equation}
\label{e:distancia} d_{\mathcal{S}}(\Phi,\Psi)=\sup_{\xi\in
\Sigma_k} d_{C^r}(\phi_\xi,\psi_\xi) + |C_\Phi-C_\Psi|,
\end{equation}
where  $\Phi=\tau\ltimes\phi_\xi$ and $\Psi=\tau\ltimes\psi_\xi$.
\end{defi}

Since we will be working with hyperbolic invariant sets generated by contracting or expanding fiber maps,
we will always assume the following hypothesis. In the case $\beta<1$ (contracting fiber maps) then $\phi_\xi(\overline{D})\subset D$ for all $\xi\in \Sigma_k$,
and if $1< \lambda $ (expanding fiber maps) then $\overline{D} \subset \phi_\xi({D})$ for all $\xi\in \Sigma_k$.

\subsection{Symbolic blender-horseshoes}
\label{ss.symbolicblenderhorseshoes}

Before introducing our main tool symbolic blender-horseshoes let us  recall the definition of a blender (with one-dimensional center contracting direction) in~\cite{BDV05,BD11}:

\begin{defi}[$cs$-blenders]
\label{d.blender}
Let $f$  be a  $C^1$-diffeomorphism of 
manifold $M$ and $\Gamma\subset M$  a transitive hyperbolic set of
$f$ with a dominated splitting of the form $E^{ss} \oplus E^{c}
\oplus E^{uu}$, where its stable bundle $E^{s}=E^{ss}\oplus E^{c}$
has dimension equal to $k\geq 2$ and $E^{c}$ is one-dimensional.
The set $\Gamma$ is a \emph{cs-blender} if  it has a
\emph{$C^1$-robust superposition region of embedding disks}, that is,
there are a
$C^1$-neighborhood $\mathcal{V}$ of $f$ and a $C^1$-open set
$\mathcal{H}$ of embeddings of $k-1$ dimensional disks into
$M$ such that for every diffeomorphism $g \in \mathcal{V}$, every
disk $H \in \mathcal{H}$ intersects the local (strong) unstable manifold
$W^{uu}_{loc}(\Gamma_g)$ of the continuation $\Gamma_g$ of $\Gamma$
for $g$.
\end{defi}

To adapt the definition of a blender to the symbolic context we need to define the family
of almost horizontal disks, which  provides the superposition region.

\begin{defi}[Almost horizontal disks]
\label{d:almost-horiz}
Given $\delta>0$, $0\leq \alpha \leq 1$ and a bounded open set $B\subset M$,
we say that a set $H^s \subset \Sigma_k \times M$ is an \emph{almost $\delta$-horizontal disk} in \mbox{$\Sigma_k\times B$}
if there are $(\zeta,z) \in \Sigma_k \times B$ and
a $(\alpha,C)$-H\"older function
$$
\text{\mbox{$h:W^s_{loc}(\zeta;\tau) \to B$} such that
$\|z-h(\xi)\|<\delta$ for all $\xi \in W^s_{loc}(\zeta;\tau)$,}$$
with
$$
C\nu^\alpha<\delta \quad \text{and} \quad
  H^s=\{(\xi,h(\xi)): \xi
\in W^s_{loc}(\zeta;\tau)\}. \vspace{0.2cm}
$$
Here
$W^s_{loc}(\zeta;\tau)\eqdef\{\xi \in \Sigma_k: \xi_i=\zeta_i, \, i\geq 0\}$ is the local stable set of $\zeta$ for $\tau$.
\end{defi}

In the case of an invariant hyperbolic set of a diffeomorphism, it
is known the existence of local (strong) unstable manifolds. An
analogous result holds for the maximal invariant set $\Gamma$ in
$\Sigma_k\times \overline{D}$ of a symbolic skew-product $\Phi$ in
$\mathcal{S}^{r,\alpha}_{k,\lambda,\beta}(D)$ with $\beta<1$ and
$\alpha>0$ (more details in Section~\ref{s.skews}).
We denote by  $W^{uu}_{loc}((\xi,x);\Phi)$ these local strong
unstable sets and define the \emph{local strong unstable set} of
$\Gamma$ 
as
 $$
    W^{uu}_{loc}(\Gamma;\Phi)= \bigcup_{x\in \Gamma} W^{uu}_{loc}((\xi,x);\Phi).
 $$
We are now ready to formulate the definition of a symbolic blender.

\begin{defi}[Symbolic $cs$-blender-horseshoe]
\label{d:symbolic-blender}
Consider a symbolic skew-product map $\Phi\in\mathcal{S}^{0,\alpha}_{k,\lambda,\beta}(D)$ with $\beta<1$ and $\alpha>0$.
The maximal invariant set $\Gamma_\Phi$ of
$\Phi$ in $\Sigma_k\times \overline{D}$  is  a \emph{symbolic
$cs$-blender-horseshoe} if there are $\delta>0$, a non-empty open set $B \subset D$, and a neighborhood $\mathcal{V}$ of $\Phi$ in $\mathcal{S}^{0,\alpha}_{k,\lambda,\beta}(D)$
such that for every $\Psi\in\mathcal{V}$ and every $\delta$-horizontal disk $H^s$ in $\Sigma_k \times B$  it holds that
\begin{equation*}
\label{e.symbolicbis}
W^{uu}_{loc}(\Gamma_\Psi;\Psi)
\cap H^s \not = \emptyset,
\end{equation*}
where $\Gamma_\Psi$ is the maximal invariant set of $\Psi$ in $\Sigma_k\times\overline{D}$.
\end{defi}
The family of $\delta$-horizontal disks in $\Sigma_k \times B$
is called the
 \emph{superposition region} of
the symbolic blenders and the open set $B$ is the \emph{superposition domain}.


A special case  of symbolic skew-product maps are the {\emph{one-step}} ones
where the fiber maps $\phi_\xi$ only depend on the coordinate $\xi_0$ of
$\xi=(\xi_i)_{i\in \mathbb{Z}} \in \Sigma_k$. So, we have  that $\phi_\xi =\phi_i$ if $\xi_0=i$,
and write
$\Phi= \tau \ltimes (\phi_1, \dots,\phi_k)$. Symbolic blenders in the one-step setting
were first introduced  in~\cite{NP11}.
In this paper, blender-horseshoes are obtained by small
perturbations in $\mathcal{S}^{0,\alpha}_{k,\lambda,\beta}(D)$
of one-step skew-product maps.




Next result gives conditions for the existence of blenders.

\begin{maintheorem}
\label{t.sym blender}
Consider $\Phi=\tau\ltimes(\phi_1,\ldots,\phi_k)\in
\mathcal{S}^{0,\alpha}_{k,\lambda,\beta}(D)$ with $\alpha>0$ and
$\nu^\alpha<\lambda< \beta<1$.
Assume that there exists an open set $B \subset D$
such that
$$
\overline{B} \subset \phi_1(B)\cup \dots \cup \phi_k(B).
$$
Then the maximal invariant set $\Gamma_\Phi$ of $\Phi$ in $\Sigma_k\times \overline{D}$ is a symbolic $cs$-blender-horseshoe for $\Phi$ whose superposition domain contains $B$.
\end{maintheorem}

\subsection{Standing notation}
\label{ss.standing}
Given a bi-sequence $\xi =( \ldots, \xi_{-1}; \xi_0, \xi_1, \ldots ) \in \Sigma_k$ the symbol
at the right of ``;'' is the ``$0$ coordinate''  of $\xi$.
We define the \emph{local stable/unstable set} of $\xi$ for $\tau$
as
\begin{align*}
 W^{s}_{loc}(\xi;\tau)=\{\zeta \in \Sigma_k: \zeta_i=\xi_i, \, i\geq 0\}  \quad \text{and}
 \quad  W^u_{loc}(\xi;\tau)&=\{\zeta \in \Sigma_k: \zeta_i=\xi_i, \, i\leq
 0\}.
\end{align*}
Given $\Phi=\tau\ltimes\phi_\xi$, for every $n>0$ and every
$(\xi,x)\in \Sigma_k\times\overline D$ we set
\begin{equation} \label{n.seq}
\begin{aligned}
   \phi^n_\xi(x) &\eqdef \phi_{\tau^{n-1}(\xi)}\circ\cdots\circ\phi_{\xi}(x), 
    \\
    \phi^{-n}_\xi(x) &\eqdef \phi^{-1}_{\tau^{-(n-1)}(\xi)}\circ\cdots\circ\phi^{-1}_{\xi}(x).
\end{aligned}
\end{equation}
Note that for every $n\ge 0$,
$$ \Phi^n(\xi,x)=
(\tau^n(\xi),\phi^n_\xi(x)) \quad \text{and}  \quad \Phi^{-n}(\xi,x)=
(\tau^{-n}(\xi),\phi^{-n}_{\tau^{-1}(\xi)}(x)).$$

\noindent For simplicity we
 write $\mathcal{S}$  in the place of $\mathcal{S}^{0,\alpha}_{k, \lambda,\beta}(D)$, where the set $D\subset M$
 and all the constants, $k\in\mathbb{N}$, $0<\lambda<\beta$, $\alpha\geq 0$,  are fixed.


\subsection{Organization of the paper}
In the next section we will study  symbolic skew-products in
$\mathcal{S}$,
and properties of the maximal invariant set in $\Sigma_k\times \overline{D}$.
In Section~\ref{s.existenceoneste}
we will look at symbolic blender-horseshoes in the one-step
setting and then
prove Theorem~\ref{t.sym blender}
in Section~\ref{s:manystep}. Finally, in Section~\ref{s:app},
we will prove Theorems~\ref{t.ciclo robusto} and~\ref{t.app mix}.
In Appendix \ref{app-reduction} we include the reduction to perturbation of symbolic skew-products.

\section{Symbolic skew-products}
\label{s.skews}
In this section given $\Phi=\tau\ltimes
\phi_\xi\in \mathcal{S}$ with $\alpha>0$, we will first study the
existence of strong stable and unstable laminations. Afterwards,
we will focus on skew-products with uniformly locally contracting
fiber maps and analyze the properties of the maximal invariant
set.

\subsection{Global strong stable and unstable laminations}

For this section let the domain $D$ of the fiber maps of $\Phi$
be the whole manifold $M$, and let us assume the global domination
condition $\nu^\alpha<\lambda<\beta<\nu^{-\alpha}$. Then the
familiar graph transform argument yields the local strong stable
and unstable laminations for the symbolic skew-product. To this
end, we define the local stable set
of $(\xi,x)$ for $\Phi$  as
\begin{align*}
W^s\big((\xi,x);\Phi\big) &\eqdef \big\{(\zeta,y)\in \Sigma_k
\times M: \lim_{n\to\infty} d(\Phi^{n}(\zeta,y),\Phi^{n}(\xi,x))=
0 \big\}.
\end{align*}

\begin{prop}[\cite{AV10}]
\label{p:inv-s} Let $\Phi \in
\mathcal{S}_{k,\lambda,\beta}^{0,\alpha}(M)$ be a symbolic
skew-product  satisfying the $s$-domination condition
$\nu^{\alpha}<\lambda$. Then, there exists a partition
$$
\mathcal{W}^{ss}_{\Phi}=\{W^{ss}_{loc}((\xi,x);\Phi): \, (\xi,x)
\in \Sigma_k \times M \}
$$
of $\Sigma_k \times M$ such that denoting $C=
C_{\Phi}(1-\lambda^{-1} \nu^\alpha )^{-1}\geq 0$, it holds
\begin{enumerate}
\item \label{item-kk1} every leaf $W^{ss}_{loc}((\xi,x);\Phi)$
 is the graph of an $\alpha$-H\"older function
 $$\gamma^{s}_{\xi,x,\Phi} : W^s_{loc}(\xi;\tau) \to M$$
 with $\alpha$-H\"older constant less or equal than $C$ (uniform on $\xi$ and
$x$),  and  it depends continuously on $\Phi$,
\item $W^{ss}_{loc}((\xi,x);\Phi)$
varies continuously with respect to $(\xi,x)$, i.e, the map
$$
(\xi,\xi',x) \mapsto \gamma^{s}_{\xi,x,\Phi}(\xi') \quad \text{is
continuous}
$$
where $(\xi, \xi')$ varies in the space of pairs of points in the
same local stable set for $\tau$,
\item \label{item-kk2} $\Phi(W^{ss}_{loc}((\xi,x);\Phi)) \subset
W^{ss}_{loc}(\Phi(\xi,x);\Phi)$ for all $(\xi,x) \in
\Sigma_k\times M$,
\item $W^{ss}_{loc}((\xi,x);\Phi) \subset
W^{s}((\xi,x);\Phi)$.
\end{enumerate}
\end{prop}

There is a dual statement of Proposition~\ref{p:inv-s} for the
symbolic skew-product $\Phi$ satisfying the $u$-dominating
condition $\beta<\nu^{-\alpha}$. Namely, there exists a partition
$ \mathcal{W}^{uu}_{\Phi} =\{W^{uu}_{loc}((\xi,x);\Phi): \,
(\xi,x) \in \Sigma_k \times M \}$ of $\Sigma_k\times M$ verifying
analogous properties. Each leaf $W^{uu}_{loc}((\xi,x);\Phi)$ of
this partition $\mathcal{W}^u_{\Phi}$ is said to be \emph{local
strong unstable set}
 throughout the point $(\xi,x)$. We define the \emph{global strong unstable set}
of $(\xi,x)$ as
\begin{equation*}
 W^{uu}
((\xi,x);\Phi ) \eqdef \displaystyle\bigcup_{n\geq 0} \Phi^{n}
\big( W^{uu}_{loc}(\Phi^{-n}(\xi,x);\Phi ) \big)\subset
W^u((\xi,x);\Phi).
\end{equation*}
For a $\Phi$-invariant set $\Gamma$ we define
\emph{local strong unstable set} of $\Gamma$ for $\Phi$ as
$$
W^{uu}_{loc}(\Gamma;\Phi)=\bigcup_{(\xi,x)\in\Gamma}W^{uu}_{loc}((\xi,x);\Phi).
$$
and the \emph{global unstable set of $\Gamma$} for $\Phi$ as
$$
W^{u}(\Gamma;\Phi)\eqdef\{(\xi,x)\in \Sigma_k \times M:
\lim_{n\to\infty}d(\Phi^{-n}(\xi,x),\Gamma_\Phi)= 0 \}.$$
In the
same manner, the \emph{global strong unstable
set} of $\Gamma$, denoted by $W^{uu}(\Gamma;\Phi)$,
is defined as the union of the global strong unstable sets of the
points in $\Gamma$ and it holds
$$
 W^{uu}_{loc}(\Gamma;\Phi) \subset W^{uu}(\Gamma;\Phi)
 \subset W^{u}(\Gamma;\Phi).
$$
Similarly, \emph{local/global (strong) stable sets} are defined
and satisfy analogous inclusions.

Let us call such skew-products with the global domination property
as partially hyperbolic. In part this is motivated by the
reduction of partially hyperbolic diffeomorphisms to the setting
of symbolic skew-products discussed in the
Appendix~\ref{app-reduction}.

\begin{defi}[Partially hyperbolic symbolic skew-products]
A symbolic skew-product $\Phi \in
\mathcal{S}^{r,\alpha}_{k,\lambda,\beta}(M)$  is called
\emph{partially hyperbolic} if $
\nu^\alpha<\lambda<1<\beta<\nu^{-\alpha}$.
 We will
denote by $\mathcal{PHS}_k^{r,\alpha}(M)$ the set of partially
hyperbolic symbolic $\alpha$-H\"older skew-products on
$\Sigma_k\times M$ with $C^r$-fiber maps.
\end{defi}


\subsection{Skew-products with  locally contracting fibers}
Suppose now that the fiber maps $\phi_\xi$ of the skew-product
$\Phi=\tau\ltimes \phi_\xi\in \mathcal{S}$ are locally contracting
on the open set $D\subset M$, that is $\beta<1$ and $D$ is forward
invariant by $\phi_\xi$. These are the hypothesis in the
definition of the symbolic blender. We will analyze the dynamics
of the maximal invariant set $\Gamma_\Phi$ in $\Sigma_k\times
\overline{D}$.

We will also need a local version of Proposition~\ref{p:inv-s}
from the previous section on the existence of the unstable
partition when there is a global $u$-domination. The same
principles used in that proof can be applied to show
the existence of a partition $\mathcal{W}^{uu}_{\Gamma_\Phi}$  of
$\Gamma_\Phi$ such that every leaf $W^{uu}_{loc}((\xi,x);\Phi)$,
$(\xi,x)\in\Gamma_\Phi$, is the graph of an $\alpha$-H\"older map
$\gamma^{u}_{\xi,x,\Phi} : W^u_{loc}(\xi;\tau) \to M$, and
$$
W^{uu}_{loc}((\xi,x);\Phi) \subset W^u((\xi,x);\Phi).
$$
In fact, $\Gamma_\Phi$ is an invariant attracting graph for $\Phi$
and, thus, $\Phi|_{\Gamma_\Phi}$ is conjugate to
$\tau|_{\Sigma_k}$ and the projection of  $\Gamma_\Phi$ in $M$
depends continuously in the Hausdorff metric on $\Phi$. These are
a version of well-known results from~\cite{HPS77} in our setting
of symbolic skew-products. We summarize these properties of the
maximal invariant set in the following theorem and for the
completeness of the paper, we give, in our symbolic setting, its
non-trivial proof  in Appendix~\ref{app-thmB}. Before stating, we
denote by $\mathrm{Per}(\Phi)$ the set of periodic points of
$\Phi$ and that $\mathscr{P}: \Sigma_k \times M \to M$ is the
projection on the fiber space.

\begin{thm}
\label{t.B} Consider $\Phi \in  \mathcal{S}^{0,\alpha}_{k,
\lambda,\beta}(D)$ with $\beta<1$ and $\alpha>0$. Then
$$
\Gamma_\Phi \eqdef \bigcap_{n\in\mathbb{Z}} \Phi^n \big(  \Sigma_k
\times \overline{D} \big)= \bigcap_{n\in\mathbb{N}} \Phi^n \big(
\Sigma_k \times \overline{D} \big).
$$
Moreover, the following holds:
\begin{enumerate}
 \item
the restriction of $\Phi$ to $\Gamma_\Phi$ is conjugate to the
full shift $\tau$ of $k$ symbols,
\item
$W^{uu}((\xi,x);\Phi)=W^u \big((\xi,x);\Phi \big)\subset
\Gamma_\Phi$ for all $(\xi,x)\in \Gamma_\Phi$,
\item
there exists a unique continuous function $g_\Phi:\Sigma_k \to
\overline{D}$ such that for every periodic point $(\vartheta,p)$
of $\Phi$ it holds,
$$
\mathscr{P}(\Gamma_\Phi)=\overline{\mathrm{Per}(\Phi)\cap D}
\eqdef K_\Phi  \in \mathcal{K}(D)
$$
and
\begin{align*}
  \qquad \quad \Gamma_\Phi = W^{uu}_{loc}(\Gamma_\Phi;\Phi) = \overline{W^{u}((\vartheta,p);\Phi))}
   =\{(\xi,g_\Phi(\xi)): \xi \in \Sigma_k\},
\end{align*}
\item \label{tB-item5}
the map $ \mathscr{L}:\mathcal{S}_{k,\lambda,\beta}^{0,\alpha}(D)
\to \mathcal{K}(\overline D)$ given by $\mathscr{L}(\Phi)=K_\Phi$
is continuous.
\end{enumerate}
\end{thm}


\section{Symbolic blenders in the one-step setting}
\label{s.existenceoneste}

Let us introduce the definition given in~\cite{NP11} for
one-step symbolic blenders.
First consider the subset  $\mathcal{Q}=\mathcal{Q}_{k,\lambda,\beta}^0(D)$
of $\mathcal{S}=\mathcal{S}_{k,\lambda,\beta}^{0,\alpha}(D)$ consisting
of one-step skew-product maps.


\begin{defi}[One-step symbolic $cs$-blender]
\label{d.sym blender-NP} Consider a one-step skew-product map
$\Phi\in\mathcal{Q}$ with $\beta<1$. The maximal invariant set
$\Gamma_\Phi$ of $\Phi$ in $\Sigma_k\times \overline{D}$ is a
\emph{one-step symbolic $cs$-blender} if there are a non-empty
open set $B \subset D$, a fixed point
$(\vartheta,p)\in\Sigma_k\times D$ of $\Phi$ and a neighborhood
$\mathcal{V}$ of $\Phi$ in $\mathcal{Q}$ such that for every $\Psi
\in \mathcal{V}$ it holds
\begin{equation}
\label{e.symbolic2}
W^{u} ( (\vartheta,p_{\Psi});\Psi)
\cap \big(W^s_{loc}(\xi; \tau) \times U \big) \not = \emptyset,
\end{equation}
for every $\xi \in \Sigma_k$ and every non-empty open subset $U$ in $B$. Here $(\vartheta,p_{\Psi})$ is the continuation of the fixed point $(\vartheta,p)$ for
$\Psi$.
\end{defi}

Observe that Definition~\ref{d:symbolic-blender} implies Definition~\ref{d.sym blender-NP}.
Let $(\vartheta, p)$ be a fixed point of $\Phi\in\mathcal{Q}$.
Note that for $(\xi,x)\in \Sigma_k\times B$ and for any
$\delta>0$, the set $W^s_{loc}(\xi;\tau)\times\{x\}$ is a $\delta$-almost horizontal disk.
Then by Theorem~\ref{t.B} and Definition~\ref{d:symbolic-blender}, the closure of $W^{u}((\vartheta,p_\Psi);\Psi)$ meets $W^{s}_{loc}(\xi;\tau) \times \{x\}$ for every
perturbation $\Psi$ of $\Phi$ in $\mathcal{Q}$.
Hence the set $W^{u}((\vartheta,p_\Psi);\Psi)$ meets $W^{s}_{loc}(\xi;\tau) \times U$,
for every open set $U$ with $x\in U$.

In this section, we prove the existence of symbolic blenders in
the one-step setting. We begin studying the relation between
one-step skew-products and their associated IFS. Given maps
$\phi_1, \ldots, \phi_k \colon \overline{D} \to D$, we denote by
$\IFS(\phi_1,\ldots, \phi_k)$, the set of all compositions of the
maps $\phi_1,\dots,\phi_k$ and will refer to this as the
associated \emph{iterated function system} (or IFS) of the
one-step map $\Phi=\tau\ltimes(\phi_1,\ldots,\phi_k)$. The
$\IFS(\Phi)= \IFS(\phi_1,\ldots,\phi_k)$  has the \emph{covering
property}  if there is an open set $B\subset D$ such that
$$
\overline{B} \subset \phi_1(B)\cup \dots \cup \phi_k(B).
$$
To construct symbolic one-step blenders we use the covering
property  and the Hutchinson attractor of the associated IFS.

\subsection{One-step skew-products and IFS's}

Consider  $\Phi=\tau\ltimes(\phi_1,\ldots,\phi_k)$ and let $\IFS(\Phi)$ be the associated iterated function system.
%
The orbit of a point $x\in M$ for $\IFS(\Phi)$ is the set
$$
\mathrm{Orb}_{\Phi}(x) \eqdef \{\phi(x):  \phi \in \IFS(\phi_1,\ldots,\phi_k)\}.
$$

Next proposition shows that if $(\vartheta, p)$ is a fixed point
of $\Phi$ then $\mathrm{Orb}_{\Phi}(p)$ is the projection onto the
fiber space of the strong unstable set of $(\vartheta, p)$. This
result was proved in \cite[Proposition 2.16]{NP11}.
A consequence of this proposition is that the density
property~\eqref{e.symbolic2} of the strong unstable set in
Definition~\ref{d.sym blender-NP} of one-step symbolic blenders is
reduced to the density of the orbit of  the ``fixed point'' $p$
for the associated iteration function system.

\begin{prop} 
\label{p.proj}
Consider $\Phi=\tau\ltimes(\phi_1,\ldots,\phi_k)$
and let $(\vartheta,p)$ be a fixed point of $\Phi$.
Then
$
\mathscr{P}(W^{uu}(\vartheta,p);\Phi)) =
\mathrm{Orb}_{\Phi}(p).
$
\end{prop}


A neighborhood $\mathcal{V}$ of a one-step map $\Phi$ in $\mathcal{Q}$ is a
neighborhood in the topology of $\mathcal{S}$ intersected with
$\mathcal{Q}$. Since the topology of $\mathcal{S}$ is induced by the
distance in~\eqref{e:distancia} and
noting that for every $\Psi \in \mathcal{Q}$
its  H\"older constant is $C_\Psi=0$, it follows that $\Psi=\tau
\ltimes (\psi_1,\ldots,\psi_k)$ and
$\Phi=\tau\ltimes(\phi_1,\ldots,\phi_k)$ are $\delta$-close if
$$
d_{\mathcal{Q}}(\Psi,\Phi)= \max_{i=1,\ldots,k}
d_{C^0}(\psi_i|_D, \phi_i|_D) <\delta.
$$

A periodic point $(\vartheta,p)$ of a  skew-product map
$\Phi=\tau\ltimes\phi_\xi$ is \emph{fiber-hyperbolic} for $\Phi$
if $p$ is a hyperbolic fixed point of $\phi^n_\vartheta$, where $n$ is the
period of $(\vartheta,p)$. We analogously define
\emph{fiber-attractors} and \emph{fiber-repellors}.

\begin{prop}
\label{p:equivalencia}
Consider $\Phi \in\mathcal{Q}$,
a non-empty open set $B \subset D$, and a fiber-hyperbolic fixed point
$(\vartheta,p) \in \Sigma_k \times D$ of $\Phi$.
The following properties are equivalent:
\begin{enumerate}
\item
\label{e:rep}
there is a neighborhood $\mathcal{V}$ of $\Phi$ in $\mathcal{Q}$ such that
for every~$\Psi \in \mathcal{V}$,
\begin{equation*}
W^{uu} \big( (\vartheta,p_{\Psi});\Psi \big) \cap
\big(W^s_{loc}(\xi; \tau) \times U \big) \not = \emptyset,
\end{equation*}
for every $\xi \in \Sigma_k$ and non-empty open set $U \subset B$. Here $(\vartheta,p_\Psi)$ is the continuation of $(\vartheta,p)$ for $\Psi$. \\
\item
\label{e:repri}
$ B \subset \overline{\mathrm{Orb}_{\Psi}(p_\Psi)}$ for every
$\Psi\in \mathcal{Q}$ close to $\Phi$.
\end{enumerate}
\end{prop}

\begin{proof}
From Proposition~\ref{p.proj}, for a fixed point
$(\vartheta,p_\Psi)$ of $\Psi$, we have that $ \mathscr{P} \big(
W^{uu} ((\vartheta,p_\Psi);\Psi )\big) =
\mathrm{Orb}_{\Psi}(p_\Psi)$. Therefore, Item~\eqref{e:rep}
implies Item~\eqref{e:repri}.

For the converse, take the fixed point $(\vartheta,p_\Psi)$ of
$\Psi=\tau\ltimes (\psi_1,\ldots,\psi_k)$ close to
$\Phi=\tau\ltimes (\phi_1,\ldots,\phi_k)$ and fix $U \subset B$
and $\xi \in \Sigma_k$. By Item~\eqref{e:repri}, there is
$\psi_{i_n}\circ \dots \circ \psi_{i_1} \in \IFS(\Psi)$ such that
the point $ x= \psi_{i_n}\circ\dots\circ \psi_{i_1}(p_\Psi) \in
U$. Take $ \zeta=(\ldots
\vartheta_{-1}\vartheta_{0},i_1,\ldots,i_n; \xi_0, \xi_1, \ldots),
$ and note that $ (\zeta, x) \in W^s_{loc}(\xi;\tau) \times U $.
It is enough to see that
 $ (\zeta, x) \in W^{uu} \big( (\vartheta,p_{\Psi});\Phi \big).  $
Since $(\vartheta,p_\Psi)$ is a fixed point of  $\Psi$,  by the
choice of $x$ we have that
$$
\Psi^{-n-1}(\zeta,x)= \big( (\ldots,\vartheta_{-1};\vartheta_0,i_1,\ldots,i_n,
\xi_0, \xi_1, \ldots),p_\Psi \big) \in W^u_{loc}(\vartheta;\tau)\times
\{p_\Psi\}.
$$
Since $W^u_{loc}(\vartheta;\tau)\times
\{p_\Psi\}=W^{uu}_{loc}((\vartheta,p_\Psi);\Psi)$ then
$$
(\zeta,x) \in  \Psi^{n+1} ( W^{uu}_{loc}((\vartheta,p_\Psi);\Psi))
\subset W^{uu}( (\vartheta,p_\Psi);\Psi).
$$
Hence
$$
(\zeta, x) \in
W^{uu}\big( (\vartheta,p_\Psi);\Psi \big)
\cap  \big( W^s_{loc}(\xi;\tau) \times U \big),
$$
completing the  proof of the proposition.
\end{proof}

\begin{rem}\label{r.sextaf}
 If $(\vartheta,p)$ in Proposition~\ref{p:equivalencia} is a
fiber-attracting fixed point of
$\Phi=\tau\ltimes(\phi_1,\ldots,\phi_k)$ with $B$ contained in the
attracting region of $p$ for $\mathrm{IFS}(\Phi)$, then
Item~\eqref{e:repri} is equivalent to
\begin{equation}
\label{eq:mot}
B \subset \overline{\mathrm{Orb}_{\Psi}(x)}, \quad \text{for every
$x\in B$ and every $\Psi
\in \mathcal{Q}$ close to $\Phi$.}
\end{equation}
\end{rem}

To see why this remark is true first note that
Equation~\eqref{eq:mot} immediately implies Item~\eqref{e:repri}
(just take $x=p$). For the converse take a perturbation
$\Psi=\tau\ltimes (\psi_1,\ldots,\psi_k)$ of $\Phi$ in
$\mathcal{Q}$, a non-empty open set $U$ in $B$, and $x\in B$. By
hypotheses, there is $\psi \in \mathrm{IFS}(\Psi)$ such that
$\psi(p_\Psi) \in U$. As $U$ is open there is a neighborhood $V$
of $p_\Psi$ such that $\psi(V)\subset U$. If $\Psi$ is close
enough to $\Phi$ then $B$ is  also in the attracting region of
$p_\Psi$ for $\psi_\vartheta = \psi_{i}$ where $\vartheta_0=i$.
Thus there is $n\in \mathbb{N}$ such that $\psi^n_{i}(x) \in V$
and hence $\psi \circ \psi_i^n(x) \in U$, proving ~\eqref{eq:mot}.

Motivated by~\eqref{eq:mot}, we give the following definition:

\begin{defi}[Blending regions]
\label{def:blending region}
Let $\Phi=\tau\ltimes (\phi_1,\ldots,\phi_k)\in \mathcal{Q}$.
A non-empty open set $B \subset M$ is
called a \emph{blending region} for $\Phi$
(or for the $\mathrm{IFS}(\Phi)$)
if for every $\Psi=\tau \ltimes (\psi_1,\ldots,\psi_k)$ close
to $\Phi$
$$
B \subset \overline{\mathrm{Orb}_{\Psi}(x)} \quad \text{for all} \ x \in B.
$$
\end{defi}



\begin{prop}
\label{prop:blender-like-set}
Let $\Phi=\tau\ltimes(\phi_1,\ldots,\phi_k) \in \mathcal{Q}$ and
consider  a blending region $B\subset D$ of $\Phi$.
Suppose that there are a hyperbolic fixed point $p \in D$ of
some $\phi_i$ and a map $\phi\in \IFS(\Phi)$ with $ \phi(p) \in B$.
Then the maximal invariant set of $\Phi$ in
$\Sigma_k\times \overline{D}$ is a one-step symbolic blender.
\end{prop}

\begin{proof}
By Proposition~\ref{p:equivalencia}, it is enough to see that
$B \subset \overline{\mathrm{Orb}_{\Psi}(p_\Psi)}$, for every
$\Psi=\tau\ltimes(\psi_1,\ldots,\psi_k)$ close to $\Phi$, where
$p_\Psi$ the continuation of $p$ for $\Psi$.
By hypothesis, there are $i_n,\dots, i_1$ such that
$\phi_{i_n}\circ\dots\circ\phi_{i_1}(p) \in B. $
Since $B$ is an open set, if $\Psi=\tau\ltimes(\psi_1,\ldots,\psi_k)$ is close
enough to $\Phi$ then $\psi_{i_n}\circ\dots\circ\psi_{i_1}(p_\Psi) \in B$.
As $B$ is a blending region for $\IFS(\Phi)$ it follows
that
$$
B \subset
\overline{\mathrm{Orb}_{\Psi}(\psi_{i_n}\circ\dots\circ\psi_{i_1}(p_\Psi))}
\subset \overline{\mathrm{Orb}_{\Psi}(p_\Psi)}.
$$
This  concludes the proof of the proposition.
\end{proof}

\subsection{Blending regions for contracting IFS: Hutchinson attractor}

Given an one-step map $\Phi=\tau \ltimes (\phi_1,\ldots,\phi_k)$,
we define $\mathrm{Per}({\mathrm{IFS}(\Phi)})$
as the projection of $\mathscr{P}(\mathrm{Per}(\Phi))$ in the fiber space,
that is the set of fixed points of the maps in  $\mathrm{IFS}(\Phi)$.
Associated with $\Phi$ or with $\mathrm{IFS}(\Phi)$,
the \emph{Hutchinson's operator} is defined by
\begin{equation*}
\mathcal{G}_{\Phi}\colon \mathcal{K}(\overline D)\to
\mathcal{K}(\overline D), \qquad \mathcal{G}_{\Phi}(A) \eqdef
\phi_1(A)\cup \dots \cup \phi_k(A),
\end{equation*}
where $\mathcal{K}(\overline D)$ denotes the set of compact subsets
of $\overline{D}$ and $A\in \mathcal{K}(\overline D)$.
If the maps $\phi_i$ are contractions, then the map $\mathcal{G}_{\Phi}$ is also a contraction.
This fact leads to the following result:

\begin{prop}[\cite{Wil71,Hut81}]
\label{p.hutch} Let $\Phi\in \mathcal{Q}$
with $\beta<1$.
Then there
exists a unique compact set $K_{\mathcal{G}_{\Phi}}\in
\mathcal{K}( \overline D)$ such that
\begin{equation*}
\label{e.hutch-atractor} K_{\mathcal{G}_{\Phi}}=\mathcal{G}_{\Phi}
(K_{\mathcal{G}_{\Phi}})=
\overline{\mathrm{Per}({\mathrm{IFS}(\Phi)})\cap D} = K_\Phi.
\end{equation*}
Moreover, the set $K_{\mathcal{G}_{\Phi}}$ depends continuously (in the set $\mathcal{Q}$) on
the map $\Phi$ and is the global attractor of
$\mathcal{G}_{\Phi}$, that is, for every $A\in \mathcal{K}( \overline D)$
it holds
$$ \displaystyle\lim_{m\to \infty}
d_H(\mathcal{G}_{\Phi}^m(A),K_{\mathcal{G}_{\Phi}})=0. $$
\end{prop}

We call the compact set $K_{\mathcal{G}_{\Phi}}$ {(in what follows
denoted by $K_\Phi$)} the \emph{Hutchinson attractor} of the
contracting one-step map $\Phi$
or of its associated
$\mathrm{IFS}(\Phi)$.

Let us recall that given $x\in D$, its orbit for $\IFS(\Phi)$ is
defined by
$$
\mathrm{Orb}_{\Phi}(x) = \{ \phi(x) \colon \phi \in  \mathrm{IFS}(\Phi)  \}
= \{  \mathcal{G}_\Phi^n(x) : n \geq 0 \}.
$$
By Proposition~\ref{p.hutch}, $\mathcal{G}_{\Phi}^m(x) \overset{m\to \infty}{\longrightarrow}
K_{\Phi}$ for all $x\in D$ and thus
$K_{\Phi} \subset \overline{\mathrm{Orb}_{\Phi}(x)}.$
We now have the following straightforward consequences of Proposition~\ref{p.hutch}:

\begin{cor}
\label{c.summingall}
Consider $\Phi=\tau\ltimes(\phi_1,\ldots,\phi_k)\in \mathcal{Q}$
with $\beta<1$
and let
$K_{\Phi}$ be the  Hutchinson attractor.
\begin{enumerate}
\item \label{(i)}
For every $A \in \mathcal{K}( \overline D)$ with $A \subset
{\mathcal{G}_{\Phi}}(A)$ it holds that
$$
A\subset K_\Phi \subset
\overline{\mathrm{Orb}_{\Phi}(x)}
\quad \text{for all} \ x\in \overline D.
$$
\item  \label{(ii)}
For every $p\in K_{\Phi}$ there is a sequence~$(\sigma_n)_{n\in
\mathbb{N}}\in \{1,\ldots,k \}^\mathbb{N}$~such~that
$$
\phi^{-1}_{\sigma_{n}}\circ\cdots\circ\phi^{-1}_{\sigma_{1}}(p)\in
K_{\Phi} \quad \text{for all} \ n\in\mathbb{N}.
$$
\item  \label{(iii)}
For each open set $V$ such that $V\cap K_{\Phi}\not = \emptyset$
there exist $n\in\mathbb{N}$ and $(i_1,\ldots,i_n)\in
\{1,\ldots,k\}^n$ such that
$\phi_{i_n}\circ\cdots\circ\phi_{i_1}(K_{\Phi})\subset V$.
\end{enumerate}
\end{cor}

These results and Proposition~\ref{prop:blender-like-set} imply
that the covering property generates one-step symbolic blenders
(Definition~\ref{d.sym blender-NP}):

\begin{cor}
\label{c.junto}
Consider $\Phi\in  \mathcal{Q}$ with $\beta<1$.
Let $B\subset \overline{D}$ be a non-empty
bounded open set satisfying the covering property for $\IFS(\Phi)$.
Then for every $\Psi \in \mathcal{Q}$
close enough to $\Phi$ one has that
$
\overline{B} \subset K_\Psi \subset \overline{\mathrm{Orb}_{\Psi}(x)}, 
\ \text{for all} \  x\in \overline{D}, $
where $K_\Psi$ is the Hutchinson attractor of $\Psi$. 
Consequently, the maximal invariant set of $\Phi$ is a one-step symbolic blender.
\end{cor}

\section{Symbolic blenders in the H\"older setting}
\label{s:manystep}

In this this section we will prove the following result:

\begin{thm}
\label{t.covering-intersection}
Consider
$\Phi=\tau\ltimes(\phi_1,\ldots,\phi_k) \in\mathcal{S}^{0,\alpha}_{k,\lambda,\beta}(D)$
with $\alpha>0$ and $\nu^\alpha<\lambda<1$ and let $B \subset D$
be an open set. Then,
$B$ satisfies the covering property
for $\IFS(\Phi)$
if and only if there are $\delta>0$ and a neighborhood $\mathcal{V}$ of $\Phi$ in
$\mathcal{S}^{0,\alpha}_{k,\lambda,\beta}(D)$ such that for every
$\Psi \in \mathcal{V}$ it holds
\begin{equation}
\label{eq:inter}
\Gamma^+_\Psi(B) \cap H^s \not=\emptyset \quad \text{for every $\delta$-horizontal disk $H^s$ in $\Sigma_k\times B$,}
\end{equation}
where $\Gamma_\Psi^+(B)$
is the forward maximal invariant set of $\Psi$ in
$\Sigma_k \times B$.
\end{thm}

The above theorem proves Theorem~\ref{t.sym blender}. Indeed,
if $\Phi=\tau\ltimes(\phi_1,\ldots,\phi_k)\in\mathcal{S}_{k,\lambda,\beta}^{0,\alpha}(D)$, with
$\beta<1$, then $\phi_i(\overline{D})\subset D$ for all $i=1,\ldots,k$. Hence, if $\Psi=\tau\ltimes \psi_\xi$ is close enough to $\Phi$ then $\psi_\xi(\overline{D})\subset D$ for all $\xi\in\Sigma_k$ and thus
\begin{equation}
\label{eq:implica}
\Gamma_\Psi^+(B)\eqdef \bigcap_{n\geq 0} \Psi^n(\Sigma_k \times B) \subset  \bigcap_{n\in \mathbb{Z}} \Psi^n(\Sigma_k \times \overline{D}) \eqdef\Gamma_\Psi.
\end{equation}
Since by Theorem~\ref{t.B} one has that $W^{uu}_{loc}(\Gamma_\Psi;\Psi)=\Gamma_\Psi$,  Theorem~\ref{t.covering-intersection} and Equation~\ref{eq:implica} imply the existence of symbolic blender (Theorem~\ref{t.sym blender}).

In order to proof Theorem~\ref{t.covering-intersection}, firstly we introduce some notation and preliminary results.

\subsection{Main lemma:} Given a finite word $\bar{\omega}=\omega_{-m}\ldots
\omega_{-1}\,\omega_0\,\omega_1\ldots \omega_n$,
where $m,n\ge 0$ and $\omega_i \in \{ 1, \dots, k  \}$,
we define the \emph{bi-lateral cylinder} by
\begin{equation*}
\label{n.seqs}
\mathcal{C}_{\bar{\omega}} \eqdef \{ \xi \in \Sigma_k\colon
\xi_{j}= \omega_{j}, \ -m \leq j \leq n \}.
\end{equation*}
Given  $\zeta \in \Sigma_k$ and a 
finite word $\bar{\omega} := \bar{\omega}_{-n}=\omega_{-n}\ldots
\omega_{-1}$, where $n\ge 1$ and  $\omega_i \in \{ 1, \dots, k
\}$, we define the {\emph{relative cylinder}} by
\begin{equation}
\label{n.cili}
\mathcal{C}_{\bar\omega}(\zeta) \eqdef\{\xi \in
W^s_{loc}(\zeta;\tau): \xi_{-i}=\omega_{-i}, \ \text{for
$i=1,\ldots,n$} \}.
\end{equation}



Recall that $\mathcal{S}=\mathcal{S}_{k, \lambda, \beta}^{0,
\alpha}(D)$ is the set of symbolic skew-product maps in
Definition~\ref{d.thesetS}. In what follows
$\nu^\alpha<\lambda<1$, $\alpha>0$ and there is no restriction on
$\beta$. In the next lemma we estimate the distance between the
backward orbits of a point $x$ when
iterated~by~different~maps~$\psi^{-1}_{\xi}$.

\begin{lem}
\label{l.noholdas} Consider $\Psi=\tau\ltimes \psi_\xi\in
\mathcal{S}$ with $\alpha>0$, $\nu^\alpha<\lambda<1$  a word
$\bar{\omega}=\omega_{-n}\ldots \omega_0 \ldots \omega_{n}$,  and
a point $x\in \overline{D}$ such that for every $\zeta\in
\mathcal{C}_{\bar{\omega}}$ one has that $
\psi^{-j}_{\tau^{-1}(\zeta)}(x) \in \overline{D}$ for every $1\leq
j\leq n$.  Then\vspace{-0.2cm}
\begin{align*}
\big \|\psi^{-i}_{\tau^{-1}(\xi)}(x)-
\psi^{-i}_{\tau^{-1}(\zeta)}(x)\big\| &< C_\Psi \, \nu^{-\alpha i}
\, \sum_{j=0}^{i-1} (\lambda^{-1} \nu^{\alpha})^{j} \,
d_{\Sigma_k}(\xi,\zeta)^\alpha
\end{align*}
for all  $1\leq i\leq n$ and all
 $\xi,\zeta \in \mathcal{C}_{\bar{\omega}}$.
\end{lem}

\begin{proof}
The proof is by induction.  For $i=1$, the H\"older inequality~\eqref{eq:Holder} and
$\xi,\zeta\in \mathcal{C}_{\bar{\omega}}$ imply
$$
\|\psi^{-1}_{\tau^{-1}(\xi)}(x) -
\psi^{-1}_{\tau^{-1}(\zeta)}(x)\| \leq C_\Psi \,
d_{\Sigma_k}(\tau^{-1}(\xi),\tau^{-1}(\zeta))^\alpha \leq
\nu^{-\alpha}\,d_{\Sigma_k}(\xi,\zeta)^\alpha.
$$
We argue inductively. Suppose that the lemma holds for $i-1$,
$i<n$:
\begin{equation}
\label{e.inducao} \big \|\psi^{-(i-1)}_{\tau^{-1}(\xi)}(x)-
\psi^{-(i-1)}_{\tau^{-1}(\zeta)}(x)\big\| < C_\Psi \,
\nu^{-\alpha(i-1)} \, \sum_{j=0}^{i-2} (\lambda^{-1}
\nu^{\alpha})^{j} \,d_{\Sigma_k}(\xi,\zeta)^\alpha,
\end{equation}
for every $\xi,\zeta \in \mathcal{C}_{\bar{\omega}}$.
We will see that the estimate also holds for $i$.
By the triangle
inequality,
\begin{align*}
\big\|\psi^{-i}_{\tau^{-1}(\xi)}(x) &-
\psi^{-i}_{\tau^{-1}(\zeta)}(x) \big\| \leq  \big\|
\psi^{-i}_{\tau^{-1}(\xi)}(x)- \psi^{-1}_{\tau^{-i}(\xi)} \circ
\psi^{-(i-1)}_{\tau^{-1}(\zeta)}(x)  \big\| \,\,+\\
&+  \big\|  \psi^{-1}_{\tau^{-i}(\xi)} \circ
\psi^{-(i-1)}_{\tau^{-1}(\zeta)}(x)-
\psi^{-i}_{\tau^{-1}(\zeta)}(x)  \big\|.
\end{align*}
Since the inverse of these functions expand at most
$1/\lambda$, the above equation is less than or equal to
\begin{equation*}
\frac{1}{\lambda}\, \big\| \psi^{-(i-1)}_{\tau^{-1}(\xi)}(x)-
\psi^{-(i-1)}_{\tau^{-1}(\zeta)}(x)  \big\| + \big\|
\psi^{-1}_{\tau^{-i}(\xi)}(y)- \psi^{-1}_{\tau^{-i}(\zeta)}(y)
\big\|,
\end{equation*}
where $y=\psi^{-(i-1)}_{\tau^{-1}(\zeta)}(x) \in \overline{D}$.
By the induction hypothesis~\eqref{e.inducao},
$$
\frac{1}{\lambda}\, \big\| \psi^{-(i-1)}_{\tau^{-1}(\xi)}(x)-
\psi^{-(i-1)}_{\tau^{-1}(\zeta)}(x)  \big\| \le C_\Psi \,
\lambda^{-1} \nu^{-\alpha(i-1)} \, \sum_{j=0}^{i-2} (\lambda^{-1}
\, \nu^{\alpha})^{j} \,d_{\Sigma_k}(\xi,\zeta)^\alpha.
$$
As $y\in \overline{D}$, applying the
H\"older inequality~\eqref{eq:Holder}
and since $\xi,\zeta \in \mathcal{C}_{\bar{ \omega}}$ we get
$$
\| \psi^{-1}_{\tau^{-i}(\xi)}(y)- \psi^{-1}_{\tau^{-i}(\zeta)}(y)
\| \le C_\Psi\,
\,d_{\Sigma_k}(\tau^{-i}(\xi),\tau^{-i}(\zeta))^\alpha \leq C_\Psi
\nu^{-\alpha i}\,d_{\Sigma_k}(\xi,\zeta)^\alpha.
$$
Putting together the previous inequalities we obtain
\begin{align*}
C_\Psi \, \lambda^{-1} & \nu^{-\alpha(i-1)} \, \sum_{j=0}^{i-2}
(\lambda^{-1} \, \nu^{\alpha})^{j}
\,d_{\Sigma_k}(\xi,\zeta)^\alpha +C_\Psi \, \nu^{-\alpha i}
\,d_{\Sigma_k}(\xi,\zeta)^\alpha
\\
&= C_\Psi \, \nu^{-\alpha i} \, \sum_{j=0}^{i-1} (\lambda^{-1} \,
\nu^{\alpha})^{j}\,d_{\Sigma_k}(\xi,\zeta)^\alpha,
\end{align*}
 ending the proof of the lemma.
\end{proof}


\subsection{Proof of Theorem~\ref{t.covering-intersection}}

Consider  $\Phi=\tau\ltimes(\phi_1,\ldots,\phi_k) \in
\mathcal{S}$ with $\alpha>0$, $\nu^\alpha<\lambda<1$  and an open
set $B\subset D$. Recall that we need to prove the following:

\begin{quotation}
\noindent
$B$ satisfies the covering property
for $\IFS(\Phi)$
$\Longleftrightarrow$  there are~\mbox{$\delta>0$} and a neighborhood $\mathcal{V}$ of $\Phi$ in
$\mathcal{S}$ such that
$\Gamma^+_\Psi(B) \cap H^s \not=\emptyset$ for every
$\Psi \in \mathcal{V}$ and every
$\delta$-horizontal disk $H^s$ in $\Sigma_k\times B$.
\end{quotation}


\smallskip

\noindent $\Longleftarrow \quad$
We will see that if the covering property  is not satisfied then
intersection~\eqref{eq:inter} is also not satisfied. If
 $B$ does not satisfy the covering property then
there is $x\in \overline{B}$ such that $x\not\in \phi_i(B)$ for all
$i=1,\ldots,k$.
We can assume that $x\in B$,
otherwise we can take an arbitrarily small perturbation
$\Psi=\tau \ltimes(\psi_1,\ldots,\psi_k)$ of $\Phi$ such that the covering property
in $B$ for $\IFS(\Psi)$ is not satisfied for a point in $B$.
The condition $x\not\in \phi_i(B)$ for all
$i=1,\ldots,k$ implies that $\Phi^{-1}(\xi,x) \not \in \Sigma_k \times B$ for all
$\xi\in \Sigma_k$ and hence
$$
   (\xi,x) \not \in \bigcap_{n\geq 0} \Phi^n(\Sigma_k\times B) = \Gamma^+_\Phi(B) \quad \text{for all} \ \xi\in\Sigma_k.
$$
Therefore $\Gamma^+_\Phi(B)$ does not meet the horizontal
disk $H^s= W^s_{loc}(\xi;\tau)\times \{x\}$,
and thus the intersection property~\eqref{eq:inter} is not verified.

\medskip

\noindent
$\Longrightarrow\quad$
We split the proof of the fact that the covering property implies the
intersection condition  into two steps.

\smallskip

\noindent {\bf{Choice of the
neighborhood $\mathcal{V}$ of $\Phi$.}}
Recall that given an open covering $\mathcal{C}$ of a compact set
$X$ of a metric space, there is a constant $L>0$, called
\emph{Lebesgue number} of $\mathcal{C}$, such that every subset of
$X$ with diameter less than $L$ is contained in some member of
$\mathcal{C}$.

Let $2L>0$ be a Lebesgue number of the open covering $\{
\phi_1(B), \dots, \phi_k(B)  \}$ of  $\overline{B}$. Note that
there are $C^0$-neigh\-bor\-hoods $\mathcal{U}_i$ of $\phi_i$ such
that the family
\begin{equation}\label{e.Bi}
B_i =  \mbox{int} \Big( \bigcap_{\psi \in \mathcal{U}_i} \psi(B)
\Big), \quad \  i=1,\ldots,k,
\end{equation}
is an open covering of $\overline{B}$. By shrinking the size of
the sets $\mathcal{U}_i$ we can assume that $L$ is
a Lebesgue number of this covering. We can also assume that any $\psi
\in\mathcal{U}_i$ is a $C^0$-$(\lambda,\beta)$-Lipschitz map
on $\overline{D}$ for all $i=1,\ldots,k$.

We take  a neighborhood $\mathcal{V}$ of $\Phi$ in $\mathcal{S}$
such that if $\Psi=\tau\ltimes \psi_\xi \in \mathcal{V}$ then
$\psi_\xi \in \mathcal{U}_{i}$ with $i=\xi_0$.
By~\eqref{e.Bi},
\begin{equation}
\label{rem-ch1}
   \psi^{-1}_{\tau^{-1}(\xi)}(\overline{B_{i}}) \subset B \quad
   \text{for all }
   \xi \in \Sigma_k \ \text{with $\xi_{-1}=i$}.
\end{equation}

Since  $\Phi$ is a
one-step map then $\phi_\xi=\phi_\zeta$ if $\xi_0=\zeta_0$, and hence we can take the H\"older constant $C_\Phi=0$.
The definition of the distance in
\eqref{e:distancia} implies $C_\Psi$ is close to $C_\Phi=0$.
By hypothesis $\nu^\alpha<\lambda$, and so shrinking the neighborhood
$\mathcal{V}$ we can assume that for every
$\Psi=\tau\ltimes \psi_\xi \in \mathcal{V}$,
\begin{equation}
\label{e:lebesguecontrol}
 C_\Psi \sum_{i=0}^{\infty} (\lambda^{-1} \nu^{\alpha})^{i}<L/2.
\end{equation}
This completes the choice of the neighborhood $\mathcal{V}$ of $\Phi$.

%

\smallskip

\noindent {\bf{Existence  of a point in $\Gamma^+_\Psi(B) \cap H^s
$.}}
Let us consider $\mathcal{V}$ to be the above neighborhood of $\Phi$.
The main step is the following proposition.

\begin{prop}\label{p.jairo}
 Consider $0<\delta <\lambda L/2$ and let $H^s$ be a
$\delta$-horizontal disk in $\Sigma_k\times B$ associated with
some $(\zeta,z)\in \Sigma_k \times B$ and $\alpha$-Hölder constant
$C\geq 0$. Then for every $\Psi=\tau\ltimes \psi_\xi \in
\mathcal{V}$ there are an infinite word
\mbox{$\bar{\omega}=\ldots\omega_{-n}\ldots\omega_{-1}$} with
$\omega_{-n} \in \{1,\ldots,k\}$ and a sequence of nested compacts
subsets $\{V_n\}$ of $M$  such that for every $n\in \mathbb{N}$ it
holds that
$$
     V_n \subset \mathscr{P}
(H^s\cap(\mathcal{C}_{\bar{\omega}_{-n}}(\zeta)\times B)),
$$
$$
\psi_{\tau^{-1}(\xi)}^{-n}(V_n) \subset B \quad \text{and} \quad
\mathrm{diam} \big( \psi_{\tau^{-1}(\xi)}^{-n}(V_n)  \big) \leq
C(\lambda^{-1}\nu^\alpha)^n
$$
for all $\xi \in \mathcal{C}_{\bar{\omega}_{-n}}(\zeta)$ where
$\bar\omega_{-n}=\omega_{-n}\ldots\omega_{-1}$.
\end{prop}

Let us see how the implication $(\Longrightarrow)$ follows from
this  proposition. Let
$$
\{x\} =\bigcap_{n\in\mathbb{N}} V_n \subset B \quad \text{and}
\quad  \{\xi\}=\bigcap_{n\in\mathbb{N}}
\mathcal{C}_{\bar{\omega}_{-n}}(\zeta)\subset
W^s_{loc}(\zeta;\tau).
$$
Observe that $(\xi,x)\in H^s$ and
$\psi^{-n}_{\tau^{-1}(\xi)}(x)\in B$ for all $n\in \mathbb{N}$.
Thus $\Psi^{-n}(\xi,x) \in \Sigma_k \times B$ for all
$n\in\mathbb{N}$ and hence $(\xi,x)\in \Gamma_\Psi^+(B)$.
Therefore $\Gamma^+_\Psi(B) \cap H^s \ne \emptyset$.

To complete the proof of Theorem~\ref{t.covering-intersection} it remains to prove the proposition.

\medskip

\begin{proof}[Proof of Proposition~\ref{p.jairo}]
Fix $\Psi=\tau\ltimes \psi_\xi \in\mathcal{V}$ and consider the
$(\alpha,C)$-Hölder map \mbox{$h:W^s_{loc}(\zeta;\tau) \to B$}
associated with the  $\delta$-horizontal disk $H^s\subset
\Sigma_k\times B$ (see Definition~\ref{d:almost-horiz}). The
construction of the nested sequence of sets $\{V_n\}$ and the
infinite word
\mbox{$\bar\omega=\ldots\omega_{-m}\ldots\omega_{-1}$} is done
inductively.
Let
$$
V \eqdef \mathscr{P}(H^s)\subset B.
$$
Note that $\mathrm{diam}(V)\leq 2\delta<L$.
By the definition of the Lebesgue number, we have that
 $V\subset B_{i_1}$ for some $i_1\in \{1,\ldots,k\}$.
Recall the definition of the relative cylinder in \eqref{n.cili} associated with $\zeta \in \Sigma_k$
and the word $\bar{\omega}_{-1}=i_1$ and consider the set
$$
V_1\eqdef\mathscr{P} \big(
H^s\cap(\mathcal{C}_{\bar{\omega}_{-1}}(\zeta) \times V) \big).
\vspace{0.2cm}
$$
By construction $V_1 \subset V \subset B_{i_1}$ and by
\eqref{rem-ch1},
\begin{equation*}
\label{e.primeirosubset}
\psi_{\tau^{-1}(\xi)}^{-1}(V_1) \subset B
\quad \text{for all $\xi\in
\mathcal{C}_{\bar\omega_{-1}}(\zeta)$.}
\end{equation*}
\begin{claim}
\label{cl.diametro} $\mathrm{diam}(V_1)\leq \delta_1\eqdef
C\,\nu^{\alpha}$.
\end{claim}

\begin{proof}
Given $x$ and $y$ in $V_1$ there are $\xi$ and $\eta$ in
$\mathcal{C}_{\bar{\omega}_{-1}}(\zeta)$ such that $x=h(\xi)$ and $y=h(\eta)$.
Since $h$ is $(\alpha, C)$-H\"older continuous,
$$
\|x-y\| =\|h(\xi)-h(\eta)\| \leq C d_{\Sigma_k}(\xi,\eta)^\alpha
\leq C \, \nu^{\alpha}=\delta_1,
$$
proving the claim.
\end{proof}

By Claim~\ref{cl.diametro} and since  the fiber-maps $\psi_\xi$
are $(\lambda,\beta)$-Lipschitz on $\overline{D}$,
$$
\mathrm{diam} \big( \psi_{\tau^{-1}(\xi)}^{-1}(V_1)  \big) \leq
\lambda^{-1}\delta_1
\quad \mbox{for all $\xi\in \mathcal{C}_{\bar{\omega}_{-1}} (\zeta)$.}
$$
Recalling that
$C \nu^{\alpha}< \delta$  (see Definition~\ref{d:almost-horiz}) we get that
$$
\lambda^{-1}\delta_1 =  \lambda^{-1} C \nu^{\alpha} < \lambda^{-1}
\delta\leq L/2.
$$
Therefore,
$$
\mathrm{diam} \big( \psi_{\tau^{-1}(\xi)}^{-1}(V_1)  \big) \leq
\lambda^{-1}\delta_1 \leq L/2.
$$


Arguing inductively, we suppose that we have constructed a finite
word \mbox{$\bar{\omega}_{-n}:= \omega_{-n} \dots \omega_{-1}$}
(the word $\bar{\omega}_{-i}$ is obtained adding the letter
$\omega_{-i}$
 to the word  $\bar{\omega}_{-i+1}$)
and closed sets $V_n\subset V_{n-1} \subset \cdots \subset V_1$
with
\begin{equation*}
\mathrm{diam}(V_n)\leq C\nu^{n\alpha} \eqdef \delta_n
\end{equation*}
and such that for every $\xi\in
\mathcal{C}_{\bar{\omega}_{-n}}(\zeta)$,
\begin{equation}
\label{e.hipo}
\psi_{\tau^{-1}(\xi)}^{-n}(V_n) \subset B
\quad \text{and} \quad
\mathrm{diam} \big( \psi_{\tau^{-1}(\xi)}^{-n}(V_n)  \big) \leq
\lambda^{-n}\delta_n.
\end{equation}
We now construct the word $\bar{\omega}_{-(n+1)}$ and the closed set
$V_{n+1} \subset V_n$ satisfying analogous inclusions and inequalities.
By \eqref{e.hipo} we have that
$$
A_n\eqdef \bigcup_{\xi \,\in \, \mathcal{C}_{\bar{\omega}_{-n}}(\zeta)   }
\psi_{\tau^{-1}(\xi)}^{-n}(V_n) \subset B.
$$

\begin{claim}\label{cl.semnome}
 $\mathrm{diam}(A_n)<L$.
\end{claim}

\begin{proof}
Given $\bar{x}$ and $\bar{y}$ in $A_n$ there are $x, y \in V_n$ and
$\xi, \eta \in \mathcal{C}_{\bar{\omega}_{-n}} (\zeta)$ such that
$\bar{x}=\psi^{-n}_{\tau^{-1}(\xi)}(x)$ and $\bar{y}=\psi^{-n}_{\tau^{-1}(\eta)}(y)$.
Then
\begin{align*}
\|\bar{x}-\bar{y}\| & =
\|\psi^{-n}_{\tau^{-1}(\xi)}(x) - \psi^{-n}_{\tau^{-1}(\eta)}(y)\|\nonumber\\
& \leq
\|\psi^{-n}_{\tau^{-1}(\xi)}(x)-\psi^{-n}_{\tau^{-1}(\eta)}(x)\|
+\|\psi^{-n}_{\tau^{-1}(\eta)}(x)-\psi^{-n}_{\tau^{-1}(\eta)}(y)\|
 \nonumber\\
& \leq C_\Psi \nu^{-\alpha n} \sum_{j=0}^{n-1}
(\lambda^{-1}\nu^{\alpha})^j  \,d_{\Sigma_k}(\xi,\eta)^\alpha+
\lambda^{-n}\delta_n
\end{align*}
where the last inequality is follows from Lemma~\ref{l.noholdas}
and the induction hypothesis~\eqref{e.hipo}. Since $\xi$ and $\eta$
belong to $\mathcal{C}_{\bar{\omega}_{-n}} (\zeta)$ then
$d_{\Sigma_k}(\xi,\eta)^\alpha \leq \nu^{\alpha n}$. Hence
$$
  \|\bar{x}-\bar{y}\| \leq C_\Psi  \sum_{j=0}^{n-1}
(\lambda^{-1}\nu^{\alpha})^j + \lambda^{-n} \delta_n \leq L/2
+\lambda^{-n} \delta_n
$$
where the last inequality follows from \eqref{e:lebesguecontrol}.
Note that
$$
\lambda^{-n} \, \delta_n = C\, (\lambda^{-1}\nu^{\alpha})^n \leq
C\, \lambda^{-1}\nu^{\alpha} < \lambda^{-1} \delta < L/2.
$$
Therefore for every pair of points  $\bar{x}, \bar{y}\in A_n$
we have $\|\bar{x}-\bar{y}\| <L$  and thus
$\mathrm{diam}(A_n)<L$, proving the claim.
\end{proof}

Since $L$ is a Lebesgue number of the covering $\{B_i\}_{i=1}^k$,
the claim implies there is $i_{n+1}\in\{1,\ldots,k\}$ such
that $A_n\subset B_{i_{n+1}}$. We let
$$
\bar{\omega}_{-(n+1)}=
i_{n+1}\omega_{-n} \dots \omega_{-1}
\quad \mbox{and} \quad
V_{n+1}=\mathscr{P}
\big( H^s\cap(\mathcal{C}_{\bar{\omega}_{-(n+1)}}(\zeta) \times V_n)  \big).
$$
Note that by construction $V_{n+1}\subset V_n$.

\begin{claim}
\label{cl.diametron+1} $\mathrm{diam}(V_{n+1})\leq
C\,\nu^{(n+1)\,\alpha} \eqdef\delta_{n+1}$.
\end{claim}

\begin{proof}
Given $x, y \in V_{n+1}$ there are $\xi, \eta \in
\mathcal{C}_{\bar{\omega}_{-(n+1)}} (\zeta)$ such that $x=h(\xi)$
and $y=h(\eta)$. From the $(\alpha,C)$-H\"older continuity of $h$
and since $\xi, \eta \in \mathcal{C}_{\bar{\omega}_{-(n+1)}}
(\zeta)$ we get $ \|x-y\|\leq C d_{\Sigma_k}(\xi,\eta)^\alpha \leq
C \nu^{(n+1)\alpha}$, concluding the proof.
\end{proof}

Using that $V_{n+1} \subset V_{n}$, $\mathrm{diam}(V_{n+1})\leq
\delta_{n+1}$, and equations~\eqref{rem-ch1} and~\eqref{e.hipo} we
get
\begin{equation*}\label{e.uffn+1}
   \psi_{\tau^{-1}(\xi)}^{-(n+1)}(V_{n+1}) \subset B
\quad \text{and} \quad
\mathrm{diam}  \Big(\psi_{\tau^{-1}(\xi)}^{-(n+1)}(V_{n+1})  \Big) \leq
\lambda^{-(n+1)}\delta_{n+1},
\end{equation*}
for all $\xi\in \mathcal{C}_{\bar{\omega}_{-(n+1)}} (\zeta)$.
Thus~\eqref{e.hipo} holds for $(n+1)$-step and we can
continue arguing inductively. This completes the construction of
the infinite word $\bar\omega$ and the sequence of nested sets
$\{V_n\}$ in the proposition, ending the proof.
\end{proof}

\medskip

The proof of Theorem~\ref{t.covering-intersection} is now complete.

\subsection{Embedded blender}



For applications of blenders we would need to have additional
fiber maps, that are not necessarily contracting or forward
invariant in the region $D$. Thus, we will embed a symbolic
blender for a one-step map $\Phi$ defined on $\Sigma_k
\times M$ into another one-step map $\hat\Phi$ on $\Sigma_d \times
M$ with $d\geq k$.

\begin{prop}[Embedded blender]
\label{p:Embedded blender} Let $\hat\Phi \in
\mathcal{PHS}^{0,\alpha}_{d}(M)$ be a partially hyperbolic
one-step skew-product map with $d\geq k$. Assume that the
restriction
$$
\Phi=\hat\Phi|_{\Sigma_k\times M}=\tau \ltimes
(\phi_1,\ldots,\phi_k)
$$
has a symbolic $cs$-blender $\Gamma^{cs}_\Phi \subset
\Sigma_k\times M$ whose superposition domain contains an open set $B\subset M$ satisfying that
$
 \overline{B}\subset \phi_1(B)\cup \dots \cup \phi_k(B).
$
Then there exists $\delta>0$ such that for  every small enough $\mathcal{S}^{0,\alpha}$-perturbation $\hat\Psi$ of $\hat\Phi$ it holds
$$
  W^{uu}_{loc}(\Gamma_\Psi^{cs}; \hat\Psi) \cap H^s \not=\emptyset
  \quad \text{for every $\delta$-horizontal disk $H^s$ in $\Sigma_d \times B$}
$$
where $\Gamma^{cs}_\Psi$ is the continuation of \
 $\Gamma_\Phi^{cs}$ for $\Psi=\hat\Psi|_{\Sigma_k\times M}$.
\end{prop}

\begin{proof}
Since $\Gamma_\Phi^{cs}$ is a symbolic blender for $\Phi$ whose
superposition domain contains the open set $B$ then
$\Gamma^+_\Phi(B)\subset \Gamma_\Phi^{cs}$. Recall that
$\Gamma^+_\Phi(B)$ is the forward maximal invariant set in
$\Sigma_k \times B$ for $\Phi$. We can enlarge the size of $B$ a
little, call this new set $\tilde{B}$, assuming that for every
small enough $\mathcal{S}^{0,\alpha}$-perturbation $\hat\Psi$ of
$\hat\Phi$, it holds
\begin{equation} \label{eq:tilde}
\Gamma^{+}_\Psi(\tilde{B}) \subset \Gamma^{cs}_\Psi, \quad \text{where $\Psi=\hat\Psi|_{\Sigma_k\times M}$.}
\end{equation}

Observe that for each $(\xi,x)\in \Sigma_d\times M$, since $\hat\Phi$ is a one-step map then
$$
   W^{uu}_{loc}((\xi,x);\hat\Phi)=W^{u}_{loc}(\xi;\tau)\times \{x\}.
$$
Since the strong unstable set depends continuously on $\hat\Phi$, we can assume that for every small enough $\mathcal{S}^{0,\alpha}$-perturbation $\hat\Psi$ of $\hat\Phi$, \begin{equation}
\label{eq:uu-B}
   W^{uu}_{loc}((\xi,x);\hat\Psi) \subset \Sigma_d \times \tilde{B} \quad \text{for all $(\xi,x)\in \Sigma_d\times B$}.
\end{equation}

A slight modification of the proof of Theorem~\ref{t.covering-intersection} shows that there is $\delta>0$ such that for every small enough  $\mathcal{S}^{0,\alpha}$-perturbation $\hat\Psi$ of $\hat\Phi$ we have
\begin{equation}
\label{eq:kdB}
  \Gamma^{+}_{\hat\Psi}(\Sigma_{k,d}^-\times B) \cap H^s\not=\emptyset
\end{equation}
for every $\delta$-horizontal disk $H^s$ in $\Sigma_d \times B$.
Here
$$
   \Gamma^+_{\hat\Psi}(\Sigma_{k,d}^-\times B)=\bigcap_{n\geq 0}
   \hat\Psi^n(\Sigma_{k,d}^-\times B)
$$
and $
  \Sigma_{k,d}^-=\{\xi=(\xi_i)_{i\in\mathbb{Z}}\in \Sigma_d : \xi_{i} \in \{1,\ldots,k\} \ \text{for $i<0$}
  \}.
$
%
This assertion is showed by repeating the same argument as in the
proof of Theorem~\ref{t.covering-intersection}, using the global
$s$-domination condition $\nu^\alpha <\lambda$ and the global
Hölder continuity of the fiber maps with respect to the base
point.

Consider a  $\mathcal{S}^{0,\alpha}$-perturbation $\hat\Psi=\tau\ltimes\psi_\xi$ of $\hat\Phi$ satisfying~\eqref{eq:tilde},~\eqref{eq:uu-B} and~\eqref{eq:kdB} and let $\Gamma_\Psi^{cs} \subset \Sigma_k \times M$ be the
continuation of $\Gamma_\Phi^{cs}$ for $\Psi=\hat\Psi|_{\Sigma_k
\times M}$.

\begin{claim}
$\Gamma^+_{\hat\Psi}(\Sigma_{k,d}^-\times B)
\subset W^{uu}_{loc}(\Gamma_\Psi^{cs};\hat\Psi)$.
\end{claim}
\begin{proof}
Fix a point $(\xi,x)\in
\Gamma^+_{\hat\Psi}(\Sigma_{k,d}^-\times B)$ we show
that
\begin{equation}
\label{eq:uukG}
   W^{uu}_{loc}((\xi,x);\hat\Psi)\cap (\Sigma_k \times \tilde{B}) \subset \Gamma_{\Psi}^{cs}.
\end{equation}
This assertion proves the claim since if $(\zeta,z)\in W^{uu}_{loc}((\xi,x);\hat\Psi)\cap (\Sigma_k \times \tilde{B})$ then
$
  (\xi,x)\in W^{uu}_{loc}((\zeta,z);\hat\Psi) \subset W^{uu}_{loc}(\Gamma_\Psi^{cs};\hat\Psi).
$

Let $(\zeta,z)\in W^{uu}_{loc}((\xi,x);\hat\Psi)\cap (\Sigma_k \times \tilde{B})$. In order to prove~\eqref{eq:uukG}, observe that from~\eqref{eq:tilde} it is enough to show that
$
    \hat\Psi^{-n}(\zeta,z)=\Psi^{-n}(\zeta,z)\in \Sigma_k \times \tilde{B}$ for all $n\geq 0$.
From the invariance of the local strong unstable lamination (dual statement of Item~\eqref{item-kk2} in Proposition~\ref{p:inv-s}) it follows that
\begin{equation}
\label{eq:3}
 \hat\Psi^{-n}(\zeta,z) \in \hat\Psi^{-n}\big(W^{uu}_{loc}((\xi,x);\hat\Psi)\big)
 \subset W^{uu}_{loc}(\hat\Psi^{-n}(\xi,x);\hat\Psi)
\end{equation}
for all $n\geq 0$. Since $(\xi,x) \in
\Gamma^+_{\hat\Psi}(\Sigma_{k,d}^-\times B)$ then
$\hat\Psi^{-n}(\xi,x)\in \Sigma_d\times B$ for all $n\geq
0$ and according to~\eqref{eq:uu-B} it follows that
\begin{equation}
\label{eq:4}
  W^{uu}_{loc}(\hat\Psi^{-n}(\xi,x);\hat\Psi) \subset \Sigma_d \times \tilde{B} \quad \text{for all $n\geq 0$}.
\end{equation}
Hence, since $\zeta \in \Sigma_k$,~\eqref{eq:3} and~\eqref{eq:4} imply that
$\hat\Psi^{-n}(\zeta,z) \in \Sigma_k \times \tilde{B}$ for all $n\geq 0$.
 Therefore $(\zeta,z)\in \Gamma^+_\Psi(\tilde{B})\subset\Gamma_\Psi^{cs}$ concluding~\eqref{eq:uukG} and the proof of the claim.
\end{proof}
To conclude the proposition it suffices to observe that the above
claim and~\eqref{eq:kdB} imply that
$W^{uu}_{loc}(\Gamma_\Psi^{cs}; \hat\Psi) \cap H^s \not=\emptyset$
for every $\delta$-horizontal disk $H^s$ in $\Sigma_d \times
B$. The proof is now complete.
\end{proof}

Observe that
the set $W^s_{loc}(\zeta;\tau)\times\{z\}$ is an horizontal disk and in the one-step case this set coincides to the local strong stable set of $(\zeta,z)$. For H\"older perturbations of one-step maps, Proposition~\ref{p:inv-s} implies that the local strong stable sets are almost horizontal disks. Then, we obtain the following remark:

\begin{rem}
\label{r.distinguish}
Let $\hat\Phi \in \mathcal{PHS}^{0,\alpha}_{d}(M)$ be a
one-step skew-product map as in the hypothesis of
Proposition~\ref{p:Embedded blender}. For every small enough
$\mathcal{S}^{0,\alpha}$-perturbation $\hat\Psi$ of $\hat\Phi$ and
every $(\xi,x)\in \Sigma_d\times B$ there exits
$(\zeta,z)\in \Gamma_\Psi^{cs}$ such that
$$
    W^{uu}_{loc}((\zeta,z);\hat\Psi)\cap W^{ss}_{loc}((\xi,x);\hat\Psi)\not=\emptyset.
$$
\end{rem}

\subsection{$cu$-blenders}
Let us now define a symbolic $cu$-blender-horseshoe. Firstly we need to
introduce the associated inverse symbolic skew-product for
$\Phi=\tau\ltimes\phi_\xi$. Given $\Phi=\tau\ltimes\phi_\xi \in
\mathcal{S}_{k,\lambda,\beta}^{0,\alpha}(D)$, the symbolic skew-product
$$
  \Phi^*=\tau\ltimes \phi^*_{\xi}\in \mathcal{S}_{k,\,\beta^{-1},\,\lambda^{-1}}^{0,\alpha}(D),
$$
where $\phi^*_{\xi}: M \to M$  given by
$ \phi^*_{\xi}(x)=\phi^{-1}_{\xi^*}(x)$,
is called the \emph{associated inverse skew-product for $\Phi$}. Here
$
\xi^*=(\ldots\xi_1; \xi_0, \xi_{-1},\ldots)
$
denotes the
conjugate sequence of $\xi=(\ldots\xi_{-1}; \xi_0,
\xi_{1},\ldots)$. Note that since $\tau(\xi)^* = \tau^{-1}(\xi^*)$,
then iterates of $\Phi^*$ correspond to iterates of
$\Phi^{-1}$.
This observation allows us to define symbolic
$cu$-blender-horseshoes for skew-products in
$\mathcal{S}^{0,\alpha}_{k,\lambda,\beta}(D)$ with $\lambda>1$ and
$\alpha>0$. Namely, \emph{symbolic $cu$-blender-horseshoe} for
$\Phi$ is defined as a symbolic $cs$-blender-horseshoe for
$\Phi^*$.

One then obtains analogous results to Theorem~\ref{t.covering-intersection} and Proposition~\ref{p:Embedded blender}
for the construction of $cu$-blenders using the inverse covering property, $\overline{B}\subset \phi_1^{-1}(B)\cup\dots\cup\phi_k^{-1}(B)$.

\section{Robust heterodimensional cycles and
mixing sets}
\label{s:app}

Our goal here is to generate robust heterodimensional cycles
(Section~\ref{ss.cycle}) and robust non-hyperbolic mixing sets
(Section~\ref{ss.mixing}) in the presence of symbolic blenders.
Let us comment on the strategy used in the proof. The idea is to
give conditions on the IFS (covering property) which will imply
that the associated one-step skew-product map has a symbolic cycle
or is mixing (see Definition~\ref{def:symb-cycle}). Then, one
studies the robustness of these properties under perturbations in
the symbolic setting (see Theorems~\ref{sym.cycle}
and~\ref{topmixthm}). Finally, the smooth realization of the
one-step skew-product shows that the corresponding property is
robust in the $C^1$-topology (Proposition~\ref{p:conj} in
appendix).

In this section, $N$ and $M$ denote compact Riemannian manifolds and
$\mathrm{id} \colon M \to M$ the identity map.


\subsection{Robust heterodimensional cycles: Proof of Theorem~\ref{t.ciclo robusto}}
\label{ss.cycle} 
We will first describe how to build
robust cycles in the symbolic setting and then transfer this to
dynamics of diffeomorphisms on manifolds, thus proving
Theorem~\ref{t.ciclo robusto}. 

\begin{defi}[Symbolic cycles]
\label{def:symb-cycle}
Let $\hat\Phi \in \mathcal{PHS}_d^{r,\alpha}(M)$ be a symbolic
skew-product with a pair of $\hat\Phi$-invariant sets
$\Gamma_1$ and $\Gamma_2$. We
say that $\hat\Phi$ has a \emph{symbolic cycle} associated with
$\Gamma_{1}$ and $\Gamma_{2}$ if their stable/unstable sets meet cyclically, that is, if
\begin{equation*}
\label{e:ciclo} W^s(\Gamma_1;\hat\Phi)\cap
W^u(\Gamma_2;\hat\Phi)\ne\emptyset \quad \mbox{and} \quad
W^u(\Gamma_1;\hat\Phi)\cap
W^s(\Gamma_2;\hat\Phi)\ne\emptyset.
\end{equation*}

%

Assuming that $\Gamma_1$ and $\Gamma_2$ have continuations, we say
that the symbolic cycle is \emph{$\mathcal{S}^{r,\alpha}$-robust}
if the cyclic intersection of stable/unstable sets holds for every
small enough $\mathcal{S}^{r,\alpha}$-perturbation of
$\hat{\Phi}$.
\end{defi}


The next proposition shows how to construct robust symbolic cycles using blenders that come from the covering property.

\begin{thm}[Symbolic cycles from blenders]
\label{sym.cycle}
Let $\phi_1, \ldots, \phi_k, \phi_{k+1}, \phi_{k+2}$ be
$(\gamma,\hat{\gamma}^{-1})$-Lipschitz $C^1$-diffeomorphisms on
$M$ and consider disjoint open sets $D_{cs},D_{cu} \subset M$ such
that
$$
   \Phi=\tau\ltimes(\phi_1,\ldots,\phi_k)\in
   \mathcal{S}^{1,\alpha}_{k,\lambda_{cs},\beta_{cs}}(D_{cs})
   \cap
   \mathcal{S}^{1,\alpha}_{k,\lambda_{cu},\beta_{cu}}(D_{cu})
$$
where
$$
   \nu^{\alpha}<\gamma\leq \lambda_{cs}<
   \beta_{cs}<1<\lambda_{cu}<\beta_{cu}\leq
   \hat{\gamma}^{-1}<\nu^{-\alpha}.
$$
Assume that there are open subsets $B_{cs} \subset D_{cs}$ and
$B_{cu} \subset D_{cu}$  such that the following properties hold:
\begin{itemize}
\item Covering property:
\begin{equation*}
\overline{B_{cs}}\subset \bigcup_{i=1}^k \phi_i(B_{cs}) \quad
\text{and} \quad \overline{B_{cu}}\subset \bigcup_{i=1}^k
\phi^{-1}_i(B_{cu})
\end{equation*}
\item Cyclic intersections: there exist  $x \in B_{cs}$, $y\in
{D}_{cu}$ and $m,n>0$  such that
$$
\phi_{k+1}^{n}(x) \in B_{cu} \quad \text{and} \quad
\phi_{k+2}^{m}(y) \in {D}_{cs}.
$$
\end{itemize}
Then the one-step
$\hat{\Phi}=\tau\ltimes(\phi_1,\ldots,\phi_k,\phi_{k+1},\phi_{k+2})
\in \mathcal{PHS}_{k+2}^{1,\alpha}(M)$ has  a $\mathcal{S}^{1,\alpha}$-robust symbolic
cycle associated with symbolic blenders whose
superposition domains contain $B_{cs}$ and $B_{cu}$.
\end{thm}

\begin{proof}
We split the proof of the  proposition into two lemmas. Firstly,
we need to choose the $\mathcal{S}^{1,\alpha}$-neighborhood
$\mathcal{V}$ of
$\hat\Phi=\tau\ltimes(\phi_1,\ldots,\phi_k,\phi_{k+1},\phi_{k+2})$
in which we will work. This neighborhood is taken as the intersection
of the following three sets:
\begin{itemize}
\item $\mathcal{V}_1$ is the neighborhood of $\hat{\Phi}$ given by Proposition~\ref{p:Embedded blender}.
\item $\mathcal{V}_2$ is the neighborhood of $\hat\Phi$ such that for every $\hat{\Psi}\in\mathcal{V}_2$
and $y \in D_{cu}$ (from the cyclic intersection hypothesis) it holds
$$
    \hat\Psi^m(\mathcal{C}^+\times \{y\}) \subset \Sigma_{k+2} \times D_{cs}
$$
where $\mathcal{C}^+$ is the set of the bi-sequence $\xi=(\xi_i)_{i\in\mathbb{Z}}\in \Sigma_{k+2}$  such that $\xi_i=k+2$ for $0\leq i \leq m$ and $\xi_i \in \{1,\ldots,k\}$ otherwise.
\item $\mathcal{V}_3$ is the neighborhood of $\hat{\Phi}$ taken so that the strong stable sets are $\delta$-almost horizontal disks, where $\delta>0$ is defined in the following manner. Suppose $z\in M$
with $\|z-\phi_{k+1}^{n}(x)\|<\delta$, then for every $\hat\Psi\in \mathcal{V}_3$ it holds
$$
   \hat\Psi^{-n}(\mathcal{C}^-\times \{z\}) \subset \Sigma_{k+2}\times B_{cs}
$$
where $\mathcal{C}^-$ is the set of the bi-sequences $\xi=(\xi_i)_{i\in\mathbb{Z}}\in \Sigma_{k+2}$  such that  $\xi_i=k+2$ for $-n\leq i <0$ and $\xi_i \in \{1,\ldots,k\}$ otherwise.
\end{itemize}

Fix $\hat\Psi\in \mathcal{V}$, let us to prove that
$$
W^s(\Gamma^{cs}_\Psi;\hat{\Psi})\cap W^u(\Gamma^{cu}_\Psi;\hat{\Psi})\neq\emptyset.
$$
where $\Gamma^{cs}_\Psi$, $\Gamma^{cu}_\Psi$ are the maximal
invariant sets in $\Sigma_d \times \overline{D_{cs}}$ and
$\Sigma_d\times \overline{D_{cu}}$ of $\hat\Psi$ respectively.
Note that $\Gamma^{cs}_\Psi$ and $\Gamma^{cu}_\Psi$  coincide,
respectively, with the maximal invariant set in $\Sigma_k \times
\overline{D_{cs}}$ and $\Sigma_k\times \overline{D_{cu}}$ for
$\Psi=\hat\Psi|_{\Sigma_k\times M}$. Hence, from the covering
property assumption and Theorem~\ref{t.sym blender}, theses two
$\hat\Psi$-invariant sets are symbolic blenders for
$\Psi$ whose superposition domains contain $B_{cs}$ and $B_{cu}$
respectively.

We observe that in the next lemma
we do not use the blender property but rather that $D_{cs}$
($D_{cu}$) are forward (backward) invariant open sets for
contractive (expanding) fiber maps of $\Phi$ and its
perturbations.

\begin{lem}\label{cycle.lemma}
It holds that $\mathcal{C}^+\times \{y\} \subset
W^s(\Gamma^{cs}_\Psi;\hat{\Psi})\cap
W^u(\Gamma^{cu}_\Psi;\hat{\Psi}).$
\end{lem}

\begin{proof}
To prove this lemma we need the following claim. Here
$\Sigma_{k,k+2}^{+}$ represents the bi-sequences
$\xi=(\xi_i)_{i\in\mathbb{Z}}\in\Sigma_{k+2}$ such that $\xi_i\in
\{1,\ldots, k\}$ for $i\geq 0$, and $\Sigma_{k,k+2}^{-}$ the
bi-sequences such that $\xi_i\in \{1,\ldots,k\}$ for $i< 0$.
\begin{claim}
$\Sigma_{k,k+2}^{+}\times D_{cs} \subset W^s(\Gamma^{cs}_\Psi;\hat\Psi)$ and $\Sigma_{k,k+2}^{-}\times D_{cu} \subset W^u(\Gamma^{cu}_\Psi;\hat\Psi)$.
\end{claim}
\begin{proof}[Proof of the claim]
We will prove the first inclusion and the second inclusion is analogous. From Theorem~\ref{t.B}, the set $\Gamma_\Psi^{cs}$ is an attractive graph and so
\begin{equation}
\label{eq:conv}
   \lim_{n\to\infty} d(\hat\Psi^n(\zeta,z), \Gamma_\Psi^{cs})= 0, \quad \text{for all $(\zeta,z)\in \Sigma_k \times D_{cs}$}.
\end{equation}
Thus if $(\xi,p) \in \Sigma_{k,k+2}^{+}\times D_{cs}$ then
$$
  \emptyset \not= W^{ss}_{loc}((\xi,p);\hat\Psi) \cap (\Sigma_k \times D_{cs}) \subset W^s(\Gamma_\Psi;\hat\Psi).
$$
Now, using the triangular inequality it holds
$$
  d(\hat\Psi^n(\xi,p),\Gamma_{\Psi}^{cs}) \leq d(\hat\Psi^n(\xi,p),\hat\Psi^{n}(\zeta,z)) +d(\hat\Psi^n(\zeta,z),\Gamma_{\Psi}^{cs})
$$
where $(\zeta,z)\in W^{ss}_{loc}((\xi,p);\hat\Psi) \cap (\Sigma_k \times D_{cs})$. According to Proposition~\ref{p:inv-s} and~\eqref{eq:conv} the right part of the above inequality converges to zero and therefore we conclude the claim.
\end{proof}

To obtain the lemma we only need to note that
$\mathcal{C}^+\times\{y\}\subset\Sigma_{k,d}^{-}\times D_{cu}$  and thus  by the above claim belongs to $W^{u}(\Gamma_\Psi^{cu};\hat\Psi)$.
On the other hand, using that $\hat\Psi\in \mathcal{V}$ and again the claim, it follows that
$$
    \hat\Psi^m(\mathcal{C}^+\times\{y\}) \subset \Sigma_{k,d}^{+}\times D_{cs} \subset W^s(\Gamma_\Psi^{cs}; \hat \Psi),
$$
which concludes the proof of the lemma.
\end{proof}

The next lemma shows the other intersection of the symbolic cycle, namely
$$
W^u(\Gamma^{cs}_{\Psi};\hat\Psi)\cap
W^s(\Gamma^{cu}_{\Psi};\hat\Psi)\neq\emptyset.
$$
This is achieved by the following stronger result as a consequence of the
blender properties.

\begin{lem} \label{cycle.lemma.A}
It holds that $W^{uu}(\Gamma^{cs}_{\Psi};\hat{\Psi})\cap
W^{ss}(\Gamma^{cu}_{\Psi};\hat{\Psi})\neq\emptyset$.
\end{lem}

\begin{proof}
Take the point $x\in B_{cs}$ from the cyclic intersection hypothesis satisfying
$\phi_{k+1}^{n}(x)\in B_{cu}$.
Since $B_{cu}\subset\mathscr{P}(\Gamma^{cu}_{\Psi})$, there exists a bi-sequence $\xi\in \Sigma_k$ such that
$(\xi, \phi^{n}_{k+1}(x))\in \Gamma^{cu}_{\Psi}$.

Consider  $\zeta=(\zeta_i)_{i\in\mathbb{Z}}\in W^s_{loc}(\xi;\tau)$ with $\zeta_i=k+1$ for
$-n\leq i<0$, and let
$$
z=\gamma^s_{\xi,\phi_{k+1}^{n}(x),\hat{\Psi}}(\zeta).
$$
Remember that
$\gamma^s$ is the graph map of the local strong stable set, which is itself an almost $\delta$-horizontal disk.
Thus, $\|z -\phi_{k+1}^{n}(x)\|<\delta$ and by the hypothesis on the neighborhood
$\mathcal{V}$ of $\hat{\Phi}$ we have that
$$
\hat{\Psi}^{-{n}}(\zeta, z)\in \Sigma_{k+2}\times B_{cs}.
$$

By the property of the blender (see Proposition~\ref{p:Embedded
blender} and Remark~\ref{r.distinguish}),
$
W^{ss}_{loc}(\hat{\Psi}^{-{n}}(\zeta, z);\hat\Psi)\cap W^{uu}_{loc}(\Gamma^{cs}_{\Psi};\hat{\Psi})\neq\emptyset.
$
Let $(\beta,t)$ be a point on this intersection, then it is in $W^{uu}_{loc}(\Gamma^{cs}_{\Psi};\hat{\Psi})$ and since $(\beta,t)\in W^{ss}_{loc}(\hat{\Psi}^{-{n}}(\zeta, z);\hat\Psi)$ then
$$\hat{\Psi}^{n}(\beta,t)\in W^{ss}_{loc}((\zeta, z);\hat\Psi)= W^{ss}_{loc}((\xi,\phi_{k+1}^{n_1}(x));\hat\Psi)\subset W^{ss}_{loc}(\Gamma^{cu}_{\Psi};\hat{\Psi}),$$
proving the lemma.
\end{proof}


Combining the above two lemmas the proof of the proposition is now
complete.
\end{proof}

Now we are ready to prove Theorem~\ref{t.ciclo robusto}.

\begin{proof}[Proof of Theorem~\ref{t.ciclo robusto}]
Fix positive constants, $\gamma, \hat{\gamma}, \lambda, \beta$
such that
$$
1/2 < \gamma \leq \lambda < \beta <1 \leq \hat{\gamma}^{-1}<2.$$
We consider an arc of $C^1$-diffeomorphisms
$\{h_\varepsilon\}_{\varepsilon \in [0,\varepsilon_0]}$,
$h_\varepsilon: M \to M$ such that $h_0=\mathrm{id}$ and for each
$\varepsilon>0$ the map $h_\varepsilon$ has two hyperbolic fixed
points $p=p(\varepsilon)$ and $q=q(\varepsilon)$ such that $p$ is a
sink, $q$ is a source,
$$W^s(p,h_\varepsilon) \pitchfork
W^u(q,h_\varepsilon) \not = \emptyset,
$$
 and
$d(p(\varepsilon),q(\varepsilon))\to 0$ as $\varepsilon \to 0$. We
also assume that there are pairwise disjoint neighborhoods
$D_{cs}=D_{cs}(\varepsilon)$ and
 $D_{cu}=D_{cu}(\varepsilon)$ of
$p$ and $q$ such that the restrictions
$h_{\varepsilon}:
\overline{D_{cs}} \to D_{cs}$ and $h^{-1}_{\varepsilon}:
\overline{D_{cu}} \to D_{cu}$ are $(\lambda,\beta)$-Lipschitz
maps.

The next useful lemma tells us how to obtain the covering property
from perturbations of a single map. See the proof in~\cite[Prop. 3.2]{NP11} and \cite{HN11}.

\begin{lem}
\label{l.translacao} Let $\phi_1: \overline{D} \to D$ be a
$(\lambda, \beta)$-Lipschitz map with $\nu^\alpha
<\lambda<\beta<1$. Then there are a natural number $k$, an open
neighborhood $B$ of the fixed point of $\phi_1$ and translations
(in local coordinates) $\phi_2, \dots, \phi_k$ of $\phi_1$ such
that $
    \overline{B} \subset \phi_1(B)\cup \dots \cup \phi_k(B).
$
Consequently, the one-step map
$$
\Phi=\tau \ltimes (\phi_1, \dots, \phi_k)\in
\mathcal{S}^{0,\alpha}_{k,\lambda,\beta}(D)
$$
has a symbolic $cs$-blender whose superposition domain
contains~$B$.
\end{lem}

We observe that the number $k$ of translations of $\phi_1$ depends on the dimension of $M$ and the contraction bound $\lambda$ of $\phi_1$.


Going back to the proof of the theorem, by Lemma~\ref{l.translacao} there are $k$ (depends only on
$\lambda$ and the dimension of $M$) and open sets $B_{cs} \subset
D_{cs}$, $B_{cu} \subset D_{cs}$  containing $p$ and $q$
respectively such that
$$
\overline{B_{cs}} \subset \phi^s_1(B_{cs}) \cup \dots \cup
\phi^u_k(B_{cs}) \quad \text{and} \quad  \overline{B_{cu}} \subset
\phi^u_1(B_{cs}) \cup \dots \cup \phi^u_k(B_{cs})
$$
where $\phi^s_1=h_\varepsilon$ on $\overline{D_{cs}}$,
$\phi^u_1=h^{-1}_\varepsilon$ on $\overline{D_{cu}}$ and
$\phi^s_i=\phi^s_i$, $\phi^u_i$ are, respectively, translations of
$\phi_1^s$, $\phi_1^u$ for all $i=2,\ldots,k$ (which depends on $\varepsilon$).

Set $\phi_1= h_\varepsilon$ and let $\phi_i$ be
$(\gamma,\hat{\gamma}^{-1})$-Lipschitz $C^1$-diffeomorphisms on
$M$ such that $\phi_i|_{\overline{D_{cs}}}=\phi^s_i$ and
$\phi_i|_{\overline{D_{cu}}}=(\phi^u_i)^{-1}$ for all
$i=2,\ldots,k$. Note that
\begin{equation*}
\overline{B_{cs}}\subset \bigcup_{i=1}^k \phi_i(B_{cs}) \quad
\text{and} \quad \overline{B_{cu}}\subset \bigcup_{i=1}^k
\phi^{-1}_i(B_{cu})
\end{equation*}
and
$$\Phi=\tau\ltimes(\phi_1,\ldots,\phi_k)\in
\mathcal{S}^{1,\alpha}_{k,\lambda_{cs},\beta_{cs}}(D_{cs})\cap
\mathcal{S}^{1,\alpha}_{k,\lambda_{cu},\beta_{cu}}(D_{cu})$$ with
$\lambda_{cs}=\lambda$, $\beta_{cs}=\beta$,
$\lambda_{cu}=\beta^{-1}$ and $\beta_{cu}=\lambda^{-1}$.

Now let us build a transition map between the regions $D_{cs}$ and $D_{cu}$.
Fix a small $\varepsilon>0$ and consider the map, in local coordinates,
$\varphi(x)=x+q-p$. From the choice of $p$ and $q$, it
holds $\varphi=\varphi_\varepsilon$ goes to the identity when
$\varepsilon\to 0$.

By hypothesis, $F$ has a Small horseshoe $\Lambda$. We consider a natural number $\ell$ such that $F^\ell|_\Lambda$ is conjugated to the full shift of $d$ symbols with $d\geq k+1$ and we take $R_1,\ldots, R_{d}$ rectangles in the
ambient manifold $N$ such that $\{R_1\cap\Lambda, \ldots, R_{d}\cap \Lambda\}$ is a Markov
partition for $F^\ell|_\Lambda$. Moreover, assume that for every unit vectors $v$
and $w$ in the stable and unstable direction of $\Lambda$
respectively it holds
$$ \mu  \leq \|DF^\ell(v)\| \leq \nu \quad \text{and} \quad
\mu  \leq \|DF^{-\ell}(w)\| \leq \nu
$$
for some positive constant $\mu\leq \nu \leq \nu^\alpha <1/2$
being $\alpha=\log \nu / \log \mu \in (0,1]$.

We now modify $f_0=F\times\mathrm{id}$ in $f_0^{\ell-1}(R_i\times M)$ to get a one-parameter family of diffeomorphisms $f_\varepsilon$ satisfying
\begin{align*}
f_{\varepsilon}^\ell |_{R_i\times M}&=F^\ell \times \phi_i  \quad \text{for $i=1,\dots,k$}\\
f_{\varepsilon}^\ell |_{R_i\times M}&=F^\ell \times \varphi  \quad\text{ for $i=k+1,\dots,d$.}
\end{align*}
Note that the locally constant skew-product diffeomorphism $f_{\varepsilon}^\ell$ restricted to the set $\Lambda\times M$ is conjugated to the
partially hyperbolic symbolic skew-product
$$\Psi_{\varepsilon}=\tau\ltimes(\phi_1,\ldots,\phi_k,\varphi,\ldots,\varphi)\in
\mathcal{PHS}_{d}^{1,\alpha}(M).
$$
 In fact, from
Proposition~\ref{p:conj}, for every small enough $C^1$-perturbation
$g$ of $f_\varepsilon$, the iterate $g^\ell$ has an invariant set $\Delta \cong
\Lambda\times M$ such that $g^\ell|_{\Delta}$ is conjugated with a
$\mathcal{S}^{1,\alpha}$-perturbation $\Psi_{g^\ell}$ of
$\Psi_{\varepsilon}$.

Observe that $\Psi_{\varepsilon}$
satisfies the assumptions in Proposition~\ref{sym.cycle}: the maps $\phi_i$
have the the covering properties in the regions $B_{cs}, B_{cu}$, and the transition maps
between these domains are created via the map $\varphi$ and using the hypothesis
that $W^s(p,h_\varepsilon) \pitchfork
W^u(q,h_\varepsilon) \not = \emptyset
$. Then the symbolic skew-product
$\Psi_{\varepsilon}$ has a $\mathcal{S}^{1,\alpha}$-robust
symbolic cycle associated with symbolic blenders.

Hence, according
to Proposition~\ref{p:conj}, using the conjugation,
$f_\varepsilon^\ell$ has a $C^1$-robust cycle associated with hyperbolic
sets $\Gamma^{cs}_\varepsilon$ and
$\Gamma^{cu}_\varepsilon$ for $f_\varepsilon^\ell$ that come from the symbolic
blenders. To compute the co-index of this cycle note that the (stable) indices of $\Gamma^{cs}_\varepsilon$ and $\Gamma_\varepsilon^{cu}$ are, respectively, $\mathrm{Ind}(\Lambda)+\mathrm{dim}(M)$ and $\mathrm{Ind(\Lambda)}$. Thus, the co-index of the heterodimensional cycle is equal to the dimension of $M$.

The one-parameter family may be taken so that the same conclusions of a heterodimensional cycle with the above co-index between two hyperbolic sets
can be made for the map $f_\varepsilon$.
This ends the proof of the theorem.
\end{proof}

\subsection{Robust mixing: proof of Theorem~\ref{t.app mix}}
\label{ss.mixing}

We will build robustly mixing examples of symbolic skew-products. One of the classical ways to create robustly transitive diffeomorphisms is to construct a map that robustly has a hyperbolic periodic point with dense stable and unstable manifolds. Hence, using
the inclination lemma (or $\lambda$-lemma) one concludes the the diffeomorphism is robustly topologically mixing. We will do the same here for the symbolic case in the following theorem.

\begin{thm} \label{topmixthm}
Let $\hat\Phi \in \mathcal{PHS}^{1,\alpha}_{d}(M)$ be a one-step map with $d>k$. Assume that the restriction
$$
\Phi=\hat\Phi|_{\Sigma_k\times M}=\tau \ltimes
(\phi_1,\ldots,\phi_k)
$$
has a symbolic $cs$-blender $\Gamma^{cs}_\Phi \subset
\Sigma_k\times M$ whose superposition domain contains an open set $B\subset M$
satisfying that
$
 \overline{B}\subset \phi_1(B)\cup \dots \cup \phi_k(B).
$

Assume that there exists an attracting fixed point of $\phi_1$ and a repelling fixed point of $\phi_{k+1}$, both in $B$, having respectively stable and unstable manifold $C^1$-robustly dense in $M$.

Then, there exists a fixed point $(\vartheta,p)\in \Gamma_\Phi^{cs}$ for $\hat\Phi$ such that its stable and unstable sets,
$$
W^s((\vartheta,p);\hat\Phi) \quad   \text{and} \quad  W^u((\vartheta,p);\hat\Phi),
$$
are $\mathcal{S}^{1,\alpha}$-robustly dense on $\Sigma_d\times M$.
Consequently $\hat\Phi$ is $\mathcal{S}^{1,\alpha}$-robustly
topologically mixing.
\end{thm}

The next sequence of propositions are necessarily to prove the above theorem.
The first proposition shows how to enlarge the topological dimension of the closure of the
strong unstable set of a blender via a fiber-repelling fixed point.  Compare to~\cite[Lemma 6.12]{BDV05}.

\begin{prop}[Blender activation]
\label{l.activation}
Let $\hat\Phi \in
\mathcal{PHS}^{0,\alpha}_{d}(M)$ be a one-step map in the
hypothesis of Proposition~\ref{p:Embedded blender}. Assume that
there exists  a fiber-repelling fixed point $(\upsilon,q)$ of
$\hat\Phi$  in $\Sigma_d\times B$. Then for every small
enough $\mathcal{S}^{0,\alpha}$-perturbation $\hat\Psi$ of
$\hat\Phi$ it holds that
$$
W^{u}((\upsilon,q_{\hat\Psi});\hat\Psi)\subset
\overline{W^{uu}(\Gamma_\Psi^{cs};\hat\Psi)}
$$
where $\Gamma_\Psi^{cs}$ is continuation of $\Gamma^{cs}_\Phi$ for
$\Psi=\hat\Psi|_{\Sigma_k\times M}$.
\end{prop}

\begin{proof}
Assume that $\hat\Phi$ satisfies the global domination
$\nu^\alpha<\lambda<1<\beta<\nu^{-\alpha}$. Consider a small
$\mathcal{S}^{0,\alpha}$-neighborhood $\mathcal{V}$ of $\hat\Phi$
such that for every $\hat\Psi \in \mathcal{V}$
$C\nu^\alpha<\delta$ where
$C=C_{\hat\Psi}(1-\lambda^{-1}\nu^\alpha)^{-1}$ and $\delta>0$ is
given in Proposition~\ref{p:Embedded blender}. Moreover, we can
take
 $\delta$ small enough such that
$B_{2\delta}(q)\subset B$.

Let $\hat\Psi \in \mathcal{V}$ and consider $(\xi,x) \in
W^{u}((\upsilon,q_{\hat\Psi});\hat\Psi)$. Let $V$ be an open neighborhood of
$(\xi,x)$.
In order to prove the proposition we need
to show that the global strong unstable set of $\Gamma^{cs}_\Psi$
meets $V$. Since $\Psi^{-n}(\xi,x)$ converges to
$(\upsilon,q_{\hat\Psi})\in \Sigma_d \times B$ then for
$n$ large enough we have that $B_\delta(\Psi^{-n}(\xi,x)) \subset
\Sigma_k \times B$. Now, we define
$$
   h_n: W^s_{loc}(\tau^{-n}(\xi);\tau) \to M, \quad h_n(\zeta)=
   \mathscr{P}\circ \hat\Psi^{-n}(\tau^{n}(\zeta),x)=
   \psi^{-n}_{\tau^{n-1}(\zeta)}(x).
$$
Let $H_n$ be the graph set of $h_n$ and notice that for $n$ large
$$
   H_n \subset \Psi^{-n}((W^{s}_{loc}(\xi;\tau)\times \{x\})\cap V).
$$
\begin{claim} \label{Holderdisc}
For $n$ large enough, $H_n$ is a $\delta$-horizontal disk in $\Sigma_d\times B$.
\end{claim}

\begin{proof}
Let $\zeta, \zeta' \in W^s_{loc}(\tau^{-n}(\xi);\tau)$ be two
bisequences with  $d_{\Sigma_d}(\zeta, \zeta')=\nu^l$. Applying
Lemma~\ref{l.noholdas} to $\tau^n(\zeta), \tau^n(\zeta')$ and the
point $x$, we obtain
\begin{align*}
 \|h_n(\zeta)-h_n(\zeta')\| &=
   \|\psi^{-n}_{\tau^{-1}(\tau^n(\zeta))}(x)-\psi^{-n}_{\tau^{-1}(\tau^n(\zeta'))}(x)\|
   \\
 &\leq
   C_\Psi\nu^{\alpha((n+l)-n)}\sum_{j=0}^{n-1}(\lambda^{-1}\nu^\alpha)^j
 \\  &<\nu^{\alpha l} C_{\hat\Psi}(1-\lambda^{-1}\nu^\alpha)^{-1}=Cd_{\Sigma_d}(\zeta, \zeta')^\alpha.
\end{align*}
 Hence, from the choice of $\delta$ we have proved that
 $h_n$
is $(\alpha,C)$-Hölder continuous map with $C\nu^\alpha<\delta$. On the other hand, for $n$ large enough $H_n\subset \Sigma_d\times B$ and therefore we conclude the proof of the
claim.
\end{proof}

From the above claim and Proposition~\ref{p:Embedded
blender}, for $n$ large enough,
$$W^{uu}_{loc}(\Gamma_\Psi^{cs};\hat\Psi)\cap H_n\neq\emptyset.$$
This implies that
$$
   \emptyset \not= W^{uu}(\Gamma_\Psi^{cs};\hat\Psi)\cap \Psi^n(H_n)
   \subset W^{uu}(\Gamma_\Psi^{cs};\hat\Psi)\cap (W^{s}_{loc}(\xi;\tau)\times \{x\})\cap V
$$
and completes the proof of the proposition.
\end{proof}

The following lemma shows how to construct a one-step map having fixed points with the stable or unstable set robustly dense.

\begin{lem} \label{l:densidad-su}
Let $\Phi =\tau\ltimes \phi_\xi  \in
\mathcal{PHS}^{1,\alpha}_{d}(M)$ be a one-step skew-product map with a fiber-repelling/attracting fixed point $(\upsilon,q)$  such that the unstable/stable manifold  of $q$ for $\phi_\upsilon$ is $C^1$-robustly dense in $M$. Then the unstable/stable set of $(\upsilon,q)$ for $\Phi$ is
$\mathcal{S}^{1,\alpha}$-robustly dense in $\Sigma_{d}\times M$.
\end{lem}

\begin{proof}
Consider a $\mathcal{S}^{1,\alpha}$-perturbation $\hat\Psi=\tau\ltimes\psi_\xi$ of $\hat\Phi$. We will prove that the unstable set of the fiber-repelling fixed point $(\upsilon,q_{\hat\Psi})$, continuation of $(\upsilon,p)$ for $\hat\Psi$, is  dense in $\Sigma_d\times M$. Similarly follows for the case of the fiber-attracting fixed point and its stable set.

Let $\mathcal{C}\times V$ be any basic open set in $\Sigma_d\times M$. By the density of the unstable set of $\upsilon$ for $\tau$ on $\Sigma_d$, there exists $\xi\in \mathcal{C} \cap W^u(\upsilon;\tau)$ and so for $n$ large enough $\tau^{-n}(\xi)$ is close to $\upsilon$.
For the fiber maps  $\psi_{\zeta}$ with $\zeta$ close enough to~$\upsilon$, there exists the fiber-repelling fixed point
$q_{\hat\Psi,\,\zeta}$, which is the continuation of $q_{\hat\Psi}$.

Let $\tilde V$ be the closure of a non-empty open set in $V$.
Fixing $n$, observe that for every $m\geq 1$ it holds that
$$
\psi^{-m-n}_{\tau^{-1}(\xi)}(\tilde{V})=
\psi^{-1}_{\tau^{-m}(\tau^{-n}(\xi))}\circ\dots\circ
\psi^{-1}_{\tau^{-1}(\tau^{-n}(\xi))}\circ\psi^{-n}_{\tau^{-1}(\xi)}(\tilde{V})
$$
where $\psi_{\tau^{-i}(\tau^{-n}(\xi))}$ is $C^1$-close to
$\psi_\upsilon$ and thus also to $\phi_\upsilon$. Since by
hypothesis the unstable manifold of $q$ for $\phi_\upsilon$ is
$C^1$-robustly dense in $M$, the same holds for the fixed point
$q_{\hat\Psi,\,{\tau^{-i}(\tau^{-n}(\xi))}}$ with respect to the
maps $\psi_{\tau^{-i}(\tau^{-n}(\xi))}$.

By the continuous dependence of compact pieces of the unstable
manifold with respect to the fiber diffeomorphisms, we may assume
that for all $m\geq 1$, the unstable manifold of
$q_{\hat\Psi,\,{\tau^{-m}(\tau^{-n}(\xi))}}$ intersects
$\psi^{-n}_{\tau^{-1}(\xi)}(\tilde{V})$.

Thus for each $m\geq 1$ there is the point $x_m$ which belongs to the intersection
$$
\psi^{-n}_{\tau^{-1}(\xi)}(\tilde{V})\cap
W^u(q_{\hat\Psi,\,{\tau^{-m}(\tau^{-n}(\xi))}};\psi_{\tau^{-m}(\tau^{-n}(\xi))}).
$$
Since $\tilde{V}$ is closed we may assume that the sequence
$\{x_m\}$ converges to a point $x$ in
$\psi^{-n}_{\tau^{-1}(\xi)}(\tilde{V})\subset
\psi^{-n}_{\tau^{-1}(\xi)}(V)$. By the H\"older continuity of
fiber maps on the base sequence,
$q_{\hat\Psi,\tau^{-m}(\tau^{-n}(\xi))}$ tends to $q_\Psi$, and
then it is not hard to show that
$\hat\Psi^{-m-n}(\tau^{-n}(\xi),x)$ goes to
$(\upsilon,q_{\hat\Psi})$ as $m\rightarrow\infty$. Thus
$$\hat\Psi^n(\tau^{-n}(\xi),x)\in (\mathcal{C}\times V)\cap W^u((\upsilon,q_{\hat\Psi});\hat\Psi),$$
concluding the proof of the lemma.
\end{proof}

The following proposition demonstrates how to create a symbolic $cs$-blender with a dense
strong unstable set.

\begin{prop}
\label{l.mix}
Let $\hat\Phi \in
\mathcal{PHS}^{1,\alpha}_{d}(M)$ be a one-step  in the
hypothesis of Proposition~\ref{p:Embedded blender} with $d>k$. Assume that the fiber map $\phi_{k+1}$ has a repelling fixed point in
$B$ with the unstable manifold $C^1$-robustly dense on $M$. Then for every small enough $\mathcal{S}^{1,\alpha}$-perturbation
$\hat\Psi$ of $\hat\Phi$ it holds that
$$
     \overline{W^{uu}(\Gamma_\Psi^{cs};\hat\Psi)}=\Sigma_{d}\times M
$$
where $\Gamma_\Psi^{cs}$ is the continuation of $\Gamma_\Phi^{cs}$
for $\Psi=\hat\Psi|_{\Sigma_k\times M}$.
\end{prop}

\begin{proof}
Let $q \in B$ be the repelling fixed point of $\phi_{k+1}$. Observe that
$(\upsilon,q)$ is a fiber-repelling fixed point of $\hat\Phi$ in
$\Sigma_{d}\times B$ where
$\upsilon=(\upsilon_i)_{i\in\mathbb{Z}} \in \Sigma_{d}$ is the
bi-sequence with $\upsilon_i=k+1$ for all $i\in\mathbb{Z}$.
By
Proposition~\ref{l.activation} it follows that for every small
enough $\mathcal{S}^{1,\alpha}$-perturbation $\hat\Psi$ of
$\hat\Phi$,
$$
W^{u}((\upsilon,q_{\hat\Psi});\hat\Psi)\subset
\overline{W^{uu}(\Gamma_\Psi^{cs};\hat\Psi)}
$$
By Lemma~\ref{l:densidad-su}, $W^u((\upsilon,q_{\hat\Psi});\hat\Psi)$ is dense in
$\Sigma_{d} \times M$ implying the density of the strong unstable set of $\Gamma^{cs}_{\hat\Psi}$ for $\hat\Psi$.
\end{proof}

\begin{lem} \label{l:incluir-uu}
Consider a symbolic skew-product $\hat\Phi \in
\mathcal{PHS}_{d}^{0,\alpha}(M)$ and let $\Gamma\subset
\Sigma_d\times M$ be a $\hat\Phi$-invariant set such that
$$ \Gamma
\subset \overline{W^{uu}((\vartheta,p);\hat\Phi)}
$$ where
$(\vartheta,p)$ is a fixed point of $\hat\Phi$. Then
$
   W^{uu}(\Gamma;\hat\Phi) \subset
   \overline{W^{uu}((\vartheta,p);\hat\Phi)}.
$
\end{lem}
\begin{proof}
Consider $(\xi,x) \in W^{uu}(\Gamma;\hat\Phi)$ and let $V$ be an open
neighborhood of $(\xi,x)$. In order to prove the lemma we need to
show that the strong unstable set of the fixed point
$(\vartheta,p)$ meets $V$. Since $(\xi,x)$ belongs to the global
strong unstable set of $\Gamma$ and $\Gamma$ is
$\hat\Phi$-invariant set then there are $n\in \mathbb{N}$ and
$(\zeta,z)\in \Gamma$ such that $
  \hat\Phi^{-n}(\xi,x)\in W^{uu}_{loc}((\zeta,z);\hat\Phi).
$ Hence
\begin{equation}\label{eq:V}
 \hat\Phi^{-n}(V) \cap
W^{uu}_{loc}((\zeta,z);\hat\Phi) \not= \emptyset.
\end{equation}
Since the global strong unstable set of $(\vartheta,p)$ lies
densely in $\Gamma$ then there exists a sequence
$((\xi^{(k)},x^{(k)}))_k\subset W^{uu}((\vartheta,p);\hat\Phi)$
converging to $(\zeta,z)$.

According to Proposition~\ref{p:inv-s},
the local strong unstable set $W^{uu}_{loc}((\xi,x);\hat\Phi)$
varies continuously with respect to $(\xi,x)$ and so
$W^{uu}_{loc}((\xi^{(k)},x^{(k)});\hat\Phi)$ converges to
$W^{uu}_{loc}((\zeta,z);\hat\Phi)$. Since $\hat\Phi^{-n}(V)$
is an open set in $\Sigma_d\times M$, then this convergence
and~\eqref{eq:V} imply that for $k$ large enough
$$
\hat\Phi^{-n}(V) \cap W^{uu}_{loc}((\xi^{(k)},x^{(k)});\hat\Phi)
\not= \emptyset.
$$
 Therefore,
$$
  \emptyset \not=V \cap
  \hat\Phi^n(W^{uu}_{loc}((\xi^{(k)},x^{(k)});\hat\Phi)) \subset
  V\cap W^{uu}((\vartheta,p);\hat\Phi),
$$
completing the proof of the lemma.
\end{proof}

Now, we are ready to give the proof of Theorem~\ref{topmixthm}.

%

\begin{proof}[Proof of Theorem~\ref{topmixthm}]
Let us denote by $(\vartheta,p)$ and $(\upsilon,q)$, respectively, the fiber-attracting and fiber-repelling fixed points of $\hat\Phi$ in $\Sigma_{d}\times B$ where the bi-sequences $\vartheta=(\vartheta_i)_{i\in\mathbb{Z}}$ and $\upsilon=(\upsilon_i)_{i\in\mathbb{Z}}$  are given by $\vartheta_i=1, \upsilon_i=k+1$ for all $i\in\mathbb{Z}$.
Observe that $(\vartheta,p)$ belongs to $\Gamma^{cs}_\Phi$.

Since by hypothesis the stable manifold of $p$ for $\phi_{1}$ is $C^1$-robustly dense on $M$, then we can apply Lemma~\ref{l:densidad-su} to conclude that the
stable set of $(\vartheta, p)$ for $\hat\Phi$  is $\mathcal{S}^{1,\alpha}$-robustly dense on $\Sigma_d\times M$.

Moreover, from Theorem~\ref{t.B}, (see also Proposition~\ref{p.refe_cube}), it holds that the strong unstable set lies $\mathcal{S}^{1,\alpha}$-robustly dense in $\Gamma^{cs}_\Phi$. Therefore, by Lemma~\ref{l:incluir-uu} it follows that
$$
W^{uu}(\Gamma^{cs}_\Phi;\hat\Phi) \subset
   \overline{W^{u}((\vartheta,p);\hat\Phi)}
$$
in a $\mathcal{S}^{1,\alpha}$-robust sense.

By Proposition~\ref{l.mix},  the strong unstable set of $\Gamma_\Phi^{cs}$ for $\hat\Phi$ is $\mathcal{S}^{1,\alpha}$-robustly dense on $\Sigma_{d}\times M$ and thus we obtain that the unstable manifold of $(\vartheta,p)$ also lies $\mathcal{S}^{1,\alpha}$-robustly dense in $\Sigma_{d}\times M$.

Therefore both the stable and unstable sets of $(\vartheta,p)$ are $\mathcal{S}^{1,\alpha}$-robustly dense in $\Sigma_d\times M$.

In order to conclude the proof of the theorem we need to prove that $\hat\Phi$ is $\mathcal{S}^{1,\alpha}$-robustly topologically mixing.
Let $\hat\Psi=\tau\ltimes\psi_\xi$ be a perturbation of $\hat\Phi$ and $(\vartheta,p_{\hat\Psi})$
the continuation of the fixed point $(\vartheta,p)$.
Consider a pair of open sets $\hat U$ and $\hat V$ in $\Sigma_d \times M$. By the density of the stable and unstable sets of $(\vartheta,p_{\hat\Psi})$ for $\hat\Psi$, we can take two points
$$
 (\eta,y)\in W^{s}((\vartheta,p_{\hat\Psi});\hat\Psi) \cap \hat U \quad \text{and} \quad
(\zeta,z)\in W^{u}((\vartheta,p_{\hat\Psi});\hat\Psi)\cap\hat V.
$$
Fix $\varepsilon>2\epsilon >0$ small enough and assume that
$\hat\Psi$ belongs to a small perturbative neighborhood of
$\hat\Phi$ such that $C_{\hat\Psi}(1-\lambda^{-1}\nu^\alpha)^{-1}
< \epsilon$ and
\begin{equation} \label{eq:uu-per}
    W^{uu}_{loc}((\xi,x);\hat\Psi) \subset \Sigma_d\times B_{\varepsilon}(p_{\hat\Psi}) \quad \text{for all $(\xi,x)\in \Sigma_d\times B_\epsilon(p_{\hat\Psi})$.}
\end{equation}
Since $(\zeta,z)$ belongs to the unstable set of
$(\vartheta,p_{\hat\Psi})$ and by Lemma~\ref{l.noholdas}, it
follows that for every $n$ large enough and $\xi \in
W^s_{loc}(\zeta;\tau)$ with $d_{\Sigma_d}(\xi,\zeta)\leq \nu^n$ we
have that
\begin{align*}
\|\psi^{-n}_{\tau^{-1}(\xi)}(z)- p_{\hat\Psi}\| &\leq
\|\psi^{-n}_{\tau^{-1}(\xi)}(z)-\psi^{-n}_{\tau^{-1}(\zeta)}(z)\|+
\|\psi^{-n}_{\tau^{-1}(\zeta)}(z)- p_{\hat\Psi}\|  \\
&\leq C_{\hat\Psi}(1-\lambda^{-1}\nu^\alpha)^{-1} + \epsilon <
2\epsilon <\varepsilon.
\end{align*}
By continuity, we can take an open set $V$ in $M$ and a cylinder $C$ in $\Sigma_d$ around of $\zeta$ of large enough length $n$  such that $C\times V \subset \hat V$ and
$$
   W^{s}_{loc}(\tau^{-n}(\zeta);\tau)\times B_\varepsilon(p_{\hat\Psi})\subset \hat\Psi^{-n}((W^s_{loc}(\zeta;\tau)\cap C)\times V) \subset \hat\Psi^{-n}(\hat V).
$$

On the other hand, since $(\eta,y)$ belongs to the stable set of $(\vartheta,p_{\hat\Psi})$ it follows that $\psi_\eta^m(y) \in B_{\epsilon}(p_{\hat\Psi})$ for all $m$ large enough. Hence,~(\ref{eq:uu-per}) implies that for every $m$ large enough
\begin{align*}
 \emptyset &\not= W^{uu}_{loc}(\hat\Psi(\eta,y);\hat\Psi) \cap
 (W^{s}_{loc}(\tau^{-n}(\zeta);\tau)\times B_\varepsilon(p_{\hat\Psi})) \\ &\subset \hat\Psi^{m}(W^{uu}_{loc}((\eta,y);\hat\Psi)\cap \hat U) \cap \hat\Psi^{-n}(\hat V).
\end{align*}
Therefore, for all $n$ and $m$ large enough $\hat\Psi^{n+m}(\hat U) \cap \hat V \not=\emptyset$, proving that $\hat\Psi$ is topologically mixing.
\end{proof}


Observe that diffeomorphisms that have fixed points with
robustly dense stable/unstable manifolds can be constructed
on any manifold $M$ by perturbations of time-one maps of gradient-like
vector fields. This completes our objective in building robustly mixing
examples in the set of symbolic skew-products (compare to \cite[Lemma~4.1]{H11} and the application in~\cite{HN11}).
Lastly we conjugate these
examples with partially hyperbolic sets on manifolds.

\begin{proof}[Proof of Theorem~\ref{t.app mix} ]
The construction of an arc of diffeomorphisms on $N\times M$ with
robustly topologically mixing non-hyperbolic sets is similar to
the proof of Theorem~\ref{t.ciclo robusto}. Thus we will just note the main differences.

As in Theorem~\ref{t.ciclo robusto}, we modify $f_0=F\times \mathrm{id}$ to get a one-parameter family of diffeomorphisms $f_\varepsilon$ satisfying
\begin{align*}
f_{\varepsilon}^\ell |_{R_i\times M}&=F^\ell \times \phi_i  \quad \text{for $i=1,\dots,k$}\\
f_{\varepsilon}^\ell |_{R_i\times M}&=F^\ell \times \varphi  \quad\text{ for $i=k+1,\ldots,d$}
\end{align*}
The maps $\phi_i=\phi_{i,\,\varepsilon}$ for $i=2,\ldots,k$ are obtained from Lemma~\ref{l.translacao} as small translations (on local coordinates) of a $C^1$-perturbation $\phi_1=\phi_{1,\,\varepsilon}$ of the identity map with an attracting fixed point $p$ with unstable manifold $C^1$-robustly dense in $M$.

Moreover, the  the iterated function system generated by $\phi_1,\ldots,\phi_k$ satisfies the covering property in an open set $B_{cs}$  with $p\in B_{cs}$. The map $\varphi=\varphi_\varepsilon$ is a $C^1$-perturbation of $\mathrm{id}$ with a repelling fixed point in $B_{cs}$ whose unstable manifold is $C^1$-robustly dense on $M$.

Notice that the restriction of $f_\varepsilon^\ell$ to $\Lambda\times M$ is conjugated to
$$
\Phi_{\varepsilon}=\tau\ltimes(\phi_1,\ldots,\phi_k,\varphi,\dots,\varphi)
\in \mathcal{PHS}^{1,\alpha}_{d}(M).
$$
This partially hyperbolic one-step satisfies the assumptions in Theorem~\ref{topmixthm} and thus $\Phi_{\varepsilon}$ is $\mathcal{S}^{1,\alpha}$-robustly topologically mixing. By Proposition~\ref{p:conj}, via the conjugation, it follows that $f_\varepsilon^\ell|_{\Lambda\times M}$ is $C^1$-robustly topologically mixing and thus the same holds for $f_\varepsilon|_{\Lambda\times M}$.
%

To show that $\Lambda\times M$ and its continuations $\Delta_g$ (for small $C^1$-perturbation $g$ of $f_\varepsilon$)
are non-hyperbolic sets, let us observe that the map $\Phi_{\varepsilon}$ by construction
has fiber-repelling and fiber-attracting fixed points. Therefore, via the conjugation,
the topologically mixing sets $\Delta_g$ on the manifold $N\times M$ will have periodic points of (stable) index
$\mathrm{Ind}(\Lambda)+\mathrm{dim}(M)$ and of index $\mathrm{Ind}(\Lambda)$, thus proving the non-hyperbolicity of these sets.

This completes the proof of the last theorem.
\end{proof}

\appendix
\section{Reduction to perturbations of symbolic skew-products}
\label{app-reduction} In this appendix we will explain how to go
from partially hyperbolic diffeomorphisms to symbolic
skew-products.

Let $F$ be a $C^2$-diffeomorphism of $N$ with a horseshoe
$\Lambda\subset N$ of $d$-legs which has a $DF$-invariant
splitting of the tangent bundle
$$T_\Lambda N = E^s_\Lambda \oplus
E^u_\Lambda$$ such that in a Riemannian metric on $N$ there is
$0<\mu<\nu<1$ satisfying
$$
   \mu  \leq \|D_zF(v)\| \leq \nu
   \quad \text{and} \quad \mu  \leq \|D_zF^{-1}(w)\| \leq \nu
$$
for all unit vectors $v\in E^s_z$, $w\in E^u_z$ and $z\in
\Lambda$. We assume that $E^s_\Lambda$, $E^u_\Lambda$ are trivial
product bundles. Both, the above assumption and the extra
regularity of $F$, are technical conditions necessary for the main
result in~\cite{IN10}, but it is conjectured that it holds without
them (see also~\cite{Go06,PSW11}). This result  is our main tool
to prove the next proposition (see Proposition~\ref{p:conj}).


Let $f$ be a $C^1$ skew-product diffeomorphism over $F$ defined as
\begin{equation}
\label{eq:skew0} f: N\times M \to N\times M, \qquad
f(z,x)=(F(z),\phi(z,x))
\end{equation}
where  $\phi(z,\cdot) : M \to M$  is a family of
$C^1$-diffeomorphisms such that there are positive numbers
$\gamma$  and $\hat{\gamma}$ satisfying
\begin{align}
\label{eq:central-domination} \gamma  < \|D_z\phi(z,\cdot)v\| <
\hat{\gamma}^{-1}  \quad \text{for all unit vector} \ v \in  T_zM \ \text{and} \ z\in \Lambda.
\end{align}
We will assume that the  skew-product
diffeomorphism~\eqref{eq:skew0} satisfies that
\begin{equation}
\label{eq:dom-ass}
   \nu + L_f < \nu^{-1}  \quad \text{and} \quad  \nu <\gamma<\hat{\gamma}^{-1} < \nu^{-1}
\end{equation}
where for a fixed $\delta>0$ small enough
$$
   L_f= \sup_{x\in M} \big\{ \, \frac{\|\phi^{\pm1}(z,x)-
   \phi^{\pm1}(z',x)\|}{d_N(z,z')}: \ z,z' \in \Lambda, \   0<d_N(z, z')<\delta   \big\} \geq 0.
$$
We say that $L_f$ is the \emph{local Lipschitz constant} of $f$
(or $f^{-1}$).
The inequalities~\eqref{eq:dom-ass} are called~\emph{modified
dominated splitting condition} in~\cite{IN10}. The first
inequality is clearly verified if $L_f=0$. The other inequalities
are the dominated conditions in the definition of partial
hyperbolicity.

\begin{defi} A skew-product diffeomorphism $f$ as~\eqref{eq:skew0}
is said to be
\emph{locally constant} if $L_f=0$.
\end{defi}


Partial hyperbolic locally constant skew-product diffeomorphisms
over a horseshoe satisfy the modified dominated splitting
condition~\eqref{eq:dom-ass}. Modified dominated splitting
condition is a $C^1$-open condition since the same property is
satisfied for any diffeomorphism $C^1$-close. However, a
$C^1$-close diffeomorphism $g$ of $f$ is a priori not a
skew-product. Under the assumptions~\eqref{eq:dom-ass}, from
Hirsch-Pugh-Shub theory~\cite{HPS77} or from the recent
works of ~\cite{Go06,IN10,PSW11} it follows that $g$ is topologically
conjugated to a skew-product. We will explain more about this.

\begin{prop}
\label{p:conj} Let $f$ be a $C^1$-skew-product diffeomorphism
as~\eqref{eq:skew0} satisfying~\eqref{eq:central-domination}
and~\eqref{eq:dom-ass}. Then given $\varepsilon>0$ small enough,
any $\varepsilon$-perturbation $g$ of $f$ in the $C^1$-topology
has a locally maximal invariant set $\Delta \subset N\times M$
homeomorphic to $\Lambda \times M$ such that \mbox{$g|_\Delta$ is
conjugated to a symbolic skew-product}
$$
\Psi_g \in\mathcal{PHS}_{d}^{1,\alpha}(M)
$$
with $\alpha=\log \nu / \log \mu>0$ and local H\"older constant
$C_{\Psi_g} = L_f+O(\varepsilon)$.
\end{prop}

\begin{rem}
In the special case of $f=F\times \mathrm{id}$,
from~\cite[Theorem~B]{Go06}, the above Proposition~\ref{p:conj} holds
without both technical assumptions, the extra regularity of the
base map $F$ and the trivialization of the vector bundles.
\end{rem}

\begin{proof}[Proof of Proposition~\ref{p:conj}]
For each $z\in \Lambda$ we consider the fiber $\mathcal{L}_z=
\{z\}\times M$. The collection $\mathcal{L}$ of these fibers is an
invariant lamination of $f$. In~\cite[Theorems~7.1]{HPS77} and
also in~\cite[Theorem~A]{IN10} and~\cite[Theorem~C]{PSW11} it is
showed that this lamination is $C^1$-persistent. The
$C^1$-persistence of such lamination means that for any
$C^1$-perturbation of $f$, there exists a lamination, $C^1$-close
to $\mathcal{L}$, which is preserved by the new dynamics, and such
that the dynamics induced on the space of the leaves remains the
same.

Given $\varepsilon>0$ small enough we take $g$ a
$C^1$-diffeomorphism $\varepsilon$-close to $f$ in the
$C^1$-topology. Note that,
$$
   g(z,x)=(\tilde{F}(z,x),\,\tilde{\phi}(z,x)),
   \qquad (z,x)\in N\times M
$$
where $\tilde{\phi}(z,\cdot):M\to M$ is also a
$C^1$-diffeomorphism such that
\begin{equation}
\label{eq:lambda-beta-contracion} \gamma \, \|x-x'\| <
\|\tilde{\phi}(z,x)-\tilde{\phi}(z,x')\| < \hat{\gamma}^{-1}\,
\|z-x'\|
\end{equation}
for all
 $x$, $x' \in M$ and $z$ in a neighborhood of $\Lambda$.

For each $z\in\Lambda$, the fiber $\mathcal{W}_{\sigma(z)}$,
continuation of $\mathcal{L}_z$ for $g$ is parametrized by the
graph of a {$C^1$-map} $Q(z,\cdot): M \to M\times N$. According
to~\cite[Theorem~A]{IN10} and since $g$ is $\varepsilon$-close to
$f$, we get
\begin{align}
\label{eq:lipschitz-1}
d_{N}(Q(z,x),Q(z,x'))&\leq O(\varepsilon) \, \|x-x'\|, \\
\label{eq:lipschitz-2} d_{C^0}(Q(z,\cdot),z)&\leq O(\varepsilon).
\end{align}
For $C^1$ maps, the $C^0$ norm of the first derivative is equal
the best Lipschitz constant. Hence, the above inequalities show
that $d_{C^1}(Q(z,\cdot),z) \leq O(\varepsilon)$.

Let
$$
\Delta=\bigcup_{z\in\Lambda} \mathcal{W}_{\sigma(z)} \subset M\times N.
$$
By sending $\mathcal{W}_{\sigma(z)}$ to $z$ we may define a continuous
projection $ P\colon \Delta \to N $ such that $P(\Delta)=\Lambda$
and $F|_\Lambda \circ P = P \circ g$. Moreover,
$$
 h: \Delta \to \Lambda \times M, \qquad  h(z,x)=(P(z,x),x)
$$
is an homeomorphism whose inverse
is
$ h^{-1}(z,x)=(Q(z,x),x). $ Let
$$
\tilde{g}=h\circ g|_{\Delta} \circ h^{-1}\colon \Lambda \times M
\to \Lambda \times M.
$$
Observe that
\begin{align*}
\tilde{g}(z,x)&= (P\circ g \circ h^{-1}(z,x),
\tilde{\phi} \circ h^{-1} (z,x)) \\
&= (F|_{\Lambda}\circ P \circ h^{-1}(z,x), \tilde{\phi}(Q(z,x),x))
=(F(z), \psi(z,x)).
\end{align*}
Thus, $\tilde{g}$ is a skew-product diffeomorphism defined on
$\Lambda\times M$ which is conjugated to $g|_\Delta$ by means of
the conjugation $h$. Since for each $z\in \Lambda$ the map
$\psi(z,\cdot)$ is a composition of $C^1$-maps, then it is a
$C^1$-map. In addition, we can easily check that the rate of
contraction and expansion of these maps are uniformly close to
$\hat{\gamma}^{-1}$ and $\gamma$ respectively. Indeed,
\begin{align*}
\|&\psi(z,x)-\psi(z,x')\| \leq \\
 &\leq \|\tilde{\phi}(Q(z,x),x)- \tilde{\phi}(Q(z,x),x')\|+
 \|\tilde{\phi}(Q(z,x),x')- \tilde{\phi}(Q(z,x'),x')\|.
\end{align*}
 By means of \eqref{eq:lambda-beta-contracion} and~\eqref{eq:lipschitz-1}
it follows that
\begin{align*}
\|\psi(z,x)-\psi(z,x')\| < \hat{\gamma}^{-1} \|x-x'\| +
O(\varepsilon)  \|x-x'\|
 \leq (\hat{\gamma}^{-1}+O(\varepsilon))  \|x-x'\|.
\end{align*}
In the same way one obtains that
\begin{align*}
\|\psi(z,x)-\psi(z,x')\| \geq
|\gamma- O(\varepsilon)|\, \|x-x'\|.
\end{align*}
Taking $\varepsilon>0$ small enough we can assume
$\hat{\gamma}^{-1}$ and $\gamma$ remain, respectively, the rates
of contraction and expansion for $\psi(z,\cdot)$. This calculation
shows that the derivatives of $\psi(z,\cdot)$ and $\phi(z,\cdot)$ are
$O(\varepsilon)$-close. Moreover, since
$Q(z,\cdot)$ is $C^1$-close to the constant function $x\mapsto z$, it follows that
$\psi(z,\cdot)$ and $\phi(z,\cdot)$ are $O(\varepsilon)$-close.
Consequently $d_{C^1}(\psi(z,\cdot), \phi(z,\cdot)) \leq
O(\varepsilon)$.

Although for each $z\in\Lambda$ the maps  $\psi(z,\cdot)$ are
$C^1$-diffeomorphisms, the map $\tilde{g}$ is not necessarily a
$C^1$-diffeomorphism since $h$ is not a $C^1$-conjugation.
However, according to~\cite[Theorem~A]{IN10} it
follows that both  $h$ and $h^{-1}$ are locally $\alpha$-H\"older
continuous maps with $\alpha=\log \nu / \log \mu >0$.
Thus, the function $\psi=\pi_M
\circ g \circ h^{-1}$ is locally $\alpha$-H\"older continuos
with respect
to the base points: if $d_N(z,z')<\delta$ it holds that
\begin{equation*}
\label{eq:holder-delta}
   d_{C^0}(\psi(z,\cdot), \psi(z',\cdot)) \leq L_g \,
   d_{C^0}(h^{-1}(z,\cdot),h^{-1}(z',\cdot)) \leq L_g C_{h}
   d_N(z,z')^{\alpha}
\end{equation*}
where $L_g$  and $C_h$ are the local Lipschitz and
$\alpha$-H\"older constants of $g$ and $h$ respectively.

Recall that $F|_\Lambda$ is conjugated to the Bernoulli shift
$\tau:\Sigma_d\to \Sigma_d$. Let \mbox{$\ell : \Sigma_d \to \Lambda$} be
the topological conjugation: $\tau = \ell^{-1}\circ F|_\Lambda
\circ \ell$.
Consider $\Sigma_d$ with the metric
given in~\eqref{e:metrica}, then
from~\cite[Theorem~2.3]{Go06} and~\cite{BGV03} it follows that
$\ell$ is a Lipschitz map. That is
$ d_M(\ell(\xi),\ell(\xi'))\leq L_\ell \, d_{\Sigma_d}(\xi,\xi'). $

Therefore, we obtain  that $\tilde{g}$ is conjugated to
a locally $\alpha$-H\"older symbolic skew-product map
$$
\Psi_g: \Sigma_d \times M \to \Sigma_d\times M, \quad
\Psi_g(\xi,x)=(\tau(\xi),\psi_\xi(x)),
$$
where $\psi_\xi = \psi(\ell(\xi),\cdot): M \to M$ is a
$(\gamma,\hat\gamma^{-1})$-Lipschitz $C^1$-diffeomorphism
with the local H\"older constant $C_{\Psi_g}= L_g C_h L_\ell^\alpha \geq 0$ and $\alpha=\log \nu
/ \log \mu>0$.

The same arguments work for a small perturbation $g^{-1}$ of $f^{-1}$
and therefore the inverse map $\Psi^{-1}_g$ is also a
skew-product with the same H\"older exponent $0<\alpha\leq 1$ and local
H\"older constant $C_\Psi\geq 0$.

To end the proof we need to show that $C_{\Psi_g}=L_f + O(\varepsilon)$. To do this
note that $L_g$ is the local Lipschitz constant of $g$ (or
$g^{-1}$) i.e.,
$$
   L_g = \sup \{\frac{d(g(z,x),g(z',x'))}{d((z,x),(z',x'))}:
   (z,x), (z',x') \in \Delta, \ d_N(z,'z)<\delta \}
$$
and then $L_g= L_f+O(\varepsilon)$. Also, the local H\"older
constant $C_h$ of $h$ (or $h^{-1}$) varies continuosly with $h$ that in turn depends continuously on $g$. In fact, in view
of~\eqref{eq:lipschitz-2} we get that $h$ and $h^{-1}$ are
close to the identity and so $C_h L_\ell^\alpha=1+O(\varepsilon)$. Thus, $C_{\Psi_g}=L_g C_h L^\alpha_\ell= L_f + O(\varepsilon)$, which concludes the proof of the proposition.
\end{proof}

Note that the above result also holds if $f$ is only defined locally as
a skew-product diffeomorphism over the rectangles coming from a Markov partition of
$\Lambda$. This fact is used to define locally constant
skew-product diffeomorphisms in the proof of Theorem~\ref{t.ciclo
robusto} and Theorem~\ref{t.app mix}.

\section{Proof of Theorem~\ref{t.B}}
\label{app-thmB} In this appendix we complete the details of the
proof of Theorem~\ref{t.B}. \vspace{0.5cm}

\noindent Next result claims the existence of a ``unique invariant
attracting graph'' on $\Sigma_k\times \overline{D}$ for
$\Phi\in \mathcal{S}$. Indeed these graphs depend continuously on
$\Phi$. The  theorem below is a reformulation of the results in
\cite{HPS77} (see also~\cite[Theorem 1.1]{Stark99} and
\cite[Section 6]{BHN} for Item~(iii)).

Given a  function $ g: \Sigma_k \to \overline{D}$ its {\emph{graph
map}} is defined by
 $$
\mathrm{graph}[g] \colon  \Sigma_k  \to \Sigma_k \times M,
\qquad \mathrm{graph}[g](\xi)= \big(\xi, g(\xi)  \big)
$$
and its \emph{graph set} by $ \Gamma_g \eqdef \mathrm{image}\big(
\mathrm{graph}[g]\big)= \{ (\xi, g (\xi)) \colon \xi \in
\Sigma_k\} \subset \Sigma_k \times \overline{D}. $

\begin{thm}[\cite{HPS77,Stark99,BHN}]
\label{t.grafo} Consider $\Phi=\tau\ltimes \phi_\xi \in
\mathcal{S}$ with $\beta<1$ and $\alpha\geq 0$. Then there exists
a unique bounded continuous function $ g_\Phi: \Sigma_k \to
\overline{D}$ such that
\begin{enumerate}
\item \label{ite-thm1} $\Phi \big(\xi,g_{\Phi}(\xi) \big)=\big(\tau(\xi), g_{\Phi}\circ \tau(\xi)\big)$,
for all $\xi \in \Sigma_k$, that is, $\Phi(\Gamma_{g_\Phi})=\Gamma_{g_\Phi}$,\\ 
\item  \label{ite2-thm1}
for every $(\xi,x) \in \Sigma_k \times \overline{D}$ and $n\geq
0$, it holds
$$
\|\phi^n_\xi(x)-g_\Phi(\tau^n(\xi))\| \leq \beta^n
\|g_\Phi(\xi)-x\|,
$$
\item the map $\mathcal{L}:\mathcal{S} \to C^0(\Sigma_k,\overline{D})$ given by $\mathcal{L}(\Phi)=g_\Phi$ is continuous.
\end{enumerate}
\end{thm}

For notational simplicity in what follows we just write
$\Gamma_\Phi$ in place of $\Gamma_{g_\Phi}$ to denote the unique
contracting invariant graph set. The next proposition shows that
$\Gamma_\Phi$ is the locally maximal invariant set in
$\Sigma_k\times \overline{D}$ and it is a consequence of previous
theorem. This completes the first item in Theorem~\ref{t.B}.

\begin{prop}
\label{l.maio2} Consider $\Phi=\tau\ltimes \phi_\xi \in
\mathcal{S}$ with $\beta<1$ and $\alpha\geq 0$. Then
\begin{enumerate}
\item \label{maio-item1}
the restriction $\Phi|_{\Gamma_\Phi}$ of $\Phi$ to the set
$\Gamma_\Phi$ is conjugate to $\tau$, and
\item
the graph set $\Gamma_\Phi$ is the (forward) maximal invariant set
in $\Sigma_k \times \overline{D}$
$$
\Gamma_\Phi = \bigcap_{n\in\mathbb{Z}} \Phi^n(\Sigma_k\times
\overline{D})=\bigcap_{n\in\mathbb{N}} \Phi^n(\Sigma_k\times
\overline{D}).
$$
\end{enumerate}
\end{prop}

The following proposition is a local version
of Proposition~\ref{p:inv-s}.

\begin{prop}
\label{p:u-lam} Consider $\Phi=\tau\ltimes\phi_\xi \in
\mathcal{S}$ with $\beta<1$, $\alpha>0$, and let $\Gamma_\Phi$ be
the maximal invariant set of $\Phi$ in  $\Sigma_k \times
\overline{D}$. Then, there exists a partition
$$
\mathcal{W}^{u}_{\Gamma_\Phi}=\{W^{uu}_{loc}((\xi,x);\Phi):  \,
(\xi,x) \in \Gamma_\Phi \}
$$
of $\Gamma_\Phi$ such that for every $(\xi,x)\in \Gamma_\Phi$ it
holds
\begin{enumerate}
\addtolength{\itemsep}{.05in}
\item \label{item-kkk1}
$W^{uu}_{loc}((\xi,x);\Phi)$ is the graph of an $\alpha$-H\"older
map $\gamma^{u}_{\xi} : W^u_{loc}(\xi;\tau) \to M$,
\item \label{i.gammaeg}
$\gamma^u_\xi(\xi)=g_\Phi(\xi)=x$,
\item \label{i.pontos}
$\gamma_\xi^u(\xi'')=\gamma_{\xi'}^u(\xi'')$ for all $\xi',\xi''
\in W^u_{loc}(\xi;\tau)$,
\item \label{item-kkk2}
$\phi^{-1}_{\tau^{-1}(\xi')}\circ \gamma^u_{\xi}(\xi') =
\gamma^u_{\tau^{-1}(\xi)}\circ \tau^{-1}(\xi')$ where  $\xi' \in
W^u_{loc}(\xi;\tau)$, and
\item \label{item-kkk3}
$W^{uu}_{loc}((\xi,x);\Phi) \subset W^u((\xi,x);\Phi)$.
\end{enumerate}
\end{prop}

\begin{proof}
Let $(\xi,x)\in\Gamma_\Phi$. Note that $x=g_\Phi(\xi)$ and that
$(\phi^{n}_{\tau^{-n}(\xi)})^{-1}= \phi^{-n}_{\tau^{-1}(\xi)}$, indeed:
$$
\phi^{-1}_{\tau^{-n}(\xi)}\circ\cdots\circ\phi^{-1}_{\tau^{-1}(\xi)}
\circ
\phi_{\tau^{-1}(\xi)}\circ\cdots \circ \phi_{\tau^{-n}(\xi)}
= \phi^{-n}_{\tau^{-1}(\xi)}\circ \phi^{n}_{\tau^{-n}(\xi)}= id.
$$

We define a sequence of maps $\gamma_{\xi}^{u,n}: W^u_{loc}(\xi;\tau) \to M$ by
$$
 \gamma_{\xi}^{u,n}(\xi')=\phi^{n}_{\tau^{-n}(\xi')}\circ
   (\phi^{n}_{\tau^{-n}(\xi)})^{-1}(x)=\phi^{n}_{\tau^{-n}(\xi')}\circ
   \phi^{-n}_{\tau^{-1}(\xi)}\circ g_\Phi(\xi).
$$
Note that $\{\gamma^{u,n}_\xi\}$ is a sequence in
the complete metric space ${C}^0(W^u_{loc}(\xi;\tau),M)$.

\begin{lem}
\label{l.cauchy}
The sequence  $\{\gamma^{u,n}_\xi\}$ is Cauchy sequence and so it converges to some  continuous map
$\gamma^u_{\xi}:W^u_{loc}(\xi;\tau) \to M$.
\end{lem}

\begin{proof}
Since $\beta, \nu<1$,
to prove the first part in the lemma it is enough to see that
\begin{equation}
\label{e.uffff}
\|\gamma_{\xi}^{u,n+1}(\xi')-\gamma_{\xi}^{u,n}(\xi')\|  \leq
C_\Phi (\beta\nu^\alpha)^{n+1}\, d_{\Sigma_k}(\xi,\xi')^{\alpha}.
\end{equation}
To prove this inequality first
recall notation in~\eqref{n.seq} and observe that Item~\eqref{ite-thm1} of
Theorem~\ref{t.grafo} implies that for every $n>0$
and $\xi \in \Sigma_k$
\begin{equation}
\label{e.tipo conj}
\phi_{\xi}^n \circ g_\Phi (\xi) = g_\Phi \circ \tau^n (\xi)
\quad \text{and} \quad
\phi_{\tau^{-1}(\xi)}^{-n}\circ g_\Phi (\xi) = g_\Phi \circ \tau^{-n} (\xi).
\end{equation}
Since $\phi_\xi(\overline{D}) \subset D$ for all $\xi\in\Sigma_k$, then for every $0<i\leq n$
$$
\phi^{n-i}_{\tau^{-n}(\xi')}\circ \phi^{-n}_{\tau^{-1}(\xi)} \circ g_\Phi(\xi) =
\phi^{n-i}_{\tau^{-n}(\xi')}\circ
g_\Phi \circ \tau^{-n} (\xi)
 \in \overline{D}.
$$
To prove the inequality in \eqref{e.uffff},
observe that $\|\gamma_{\xi}^{u,n+1}(\xi')-\gamma_{\xi}^{u,n}(\xi')\| $  is equal to
\begin{align*}
&
\| \phi^{n+1}_{\tau^{-n-1}(\xi')}\circ
\phi^{-n-1}_{\tau^{-1}(\xi)}\circ g_\Phi(\xi) -
\phi^{n}_{\tau^{-n}(\xi')}\circ \phi^{-n}_{\tau^{-1}(\xi)}\circ g_\Phi(\xi)  \| \\
& \quad \leq
 \beta^{n} \|\phi_{\tau^{-n-1}(\xi')}\circ
\phi^{-1}_{\tau^{-n-1}(\xi)} \circ \phi^{-n}_{\tau^{-1}(\xi)}\circ g_\Phi(\xi) \, - \\
& \qquad \qquad \qquad \qquad \qquad \qquad  \,
 \phi_{\tau^{-n-1}(\xi')} \circ \phi^{-1}_{\tau^{-n-1}(\xi')}\circ
 \phi^{-n}_{\tau^{-1}(\xi)}\circ  g_\Phi(\xi)\| \\
& \quad \leq
 \beta^{n+1} \|\phi^{-1}_{\tau^{-n-1}(\xi)}\circ
 \phi^{-n}_{\tau^{-1}(\xi)} \circ g_\Phi(\xi) -
 \phi^{-1}_{\tau^{-n-1}(\xi')}\circ
\phi^{-n}_{\tau^{-1}(\xi)}\circ g_\Phi(\xi)\|.
\end{align*}
As $\phi^{-n}_{\tau^{-1}(\xi)}\circ g_\Phi(\xi) \in \overline{D}$,
then local H\"older continuity implies
\begin{align*}
\|\gamma_{\xi}^{u,n+1}(\xi')-\gamma_{\xi}^{u,n}(\xi')\| &\leq
 \beta^{n+1} \, C_{\Phi} \, d_{\Sigma_k} \big(  \tau^{-n-1}(\xi), \tau^{-n-1}(\xi')  \big)^\alpha\\
&\leq
  C_\Phi (\beta\nu^\alpha)^{n+1}\, d_{\Sigma_k}(\xi,\xi')^{\alpha}.
\end{align*}
The proof of the lemma is now complete.
\end{proof}


To prove that $\mathcal{W}^u_{\Gamma_\Phi}$ is a partition of $\Gamma_\Phi$ we need to show that $W^{uu}_{loc}((\xi,x);\Phi)$ is contained in $\Gamma_\Phi$ for all $(\xi,x)\in\Gamma_\Phi$.
To do this let
$$
(\xi',x')=(\xi', \gamma^u_{\xi} (\xi')) \in
W^{uu}_{loc}((\xi,x);\Phi). $$
 From~\eqref{e.tipo conj} and noting
that $x=g_\Phi(\xi)$ we have that
\begin{align*}
\|g_\Phi(\xi')-\gamma^u_{\xi}(\xi')\| &= \lim_{n\to \infty}
\|\phi^n_{\tau^{-n}(\xi')}\circ g_\Phi \circ
\tau^{-n}(\xi')-\phi^n_{\tau^{-n}(\xi')}\circ g_\Phi \circ
\tau^{-n}(\xi)\|.
\end{align*}
Since the maps $\phi_\xi$ are $\beta$-contractions and  $g_\Phi:\Sigma_k\to \overline{D}$,
if $K$ is an upper bound of the diameter of $\overline{D}$
we get that
\begin{align*}
\|g_\Phi(\xi')-\gamma^u_{\xi}(\xi')\| \leq \lim_{n\to \infty}
\beta^{n} \| g_\Phi \circ \tau^{-n}(\xi')- g_\Phi \circ
\tau^{-n}(\xi)\| \leq  \lim_{n\to \infty} K\beta^n =0.
\end{align*}
Thus
$g_\Phi(\xi')=\gamma^u_{\xi}(\xi')=x'$ and so $(\xi',x')\in
\Gamma_\Phi$. That is,
$$
W^{uu}_{loc}((\xi,x);\Phi) \subset \Gamma_\Phi
\quad \text{for all} \ (\xi,x)\in\Gamma_\Phi.
$$

Indeed, we have proved that $\gamma^u_\xi(\xi')=g_\Phi(\xi')$ for all
$\xi' \in W^u_{loc}(\xi;\tau)$, proving~\eqref{i.gammaeg}.
In particular, $\gamma^u_{\xi}(\xi'')=g_\Phi(\xi'')=\gamma^u_{\xi'}(\xi'')$
for all $\xi',\xi''\in W^u_{loc}(\xi;\tau)$, which proves \eqref{i.pontos}.

To prove that $\gamma_{\xi}^{u}$ is an $\alpha$-H\"older map (Item \eqref{item-kkk1})
first note that  by the triangle inequality
\begin{equation}
\label{eq:close-id0}
   \|\gamma_{\xi}^{u,n}(\xi')-x\| =\|\gamma_{\xi}^{u,n}(\xi')-\gamma_{\xi}^{u,n}(\xi)\| \leq \sum_{i=1}^{n} s_i(\xi')
\end{equation}
where $x=g_\Phi(\xi)$ and  $s_i(\xi')$ is given by
\begin{align*}
&\|\phi^{n-i}_{\tau^{-n+i}(\xi')}\circ \phi_{\tau^{-n-1+i}(\xi')}
\circ \phi_{\tau^{-n}(\xi)}^{i-1}\circ \phi_{\tau^{-1}(\xi)}^{-n}(x) \,- \\
&\hspace*{4cm} \phi^{n-i}_{\tau^{-n+i}(\xi')}\circ \phi_{\tau^{-n-1+i}(\xi)}
\circ \phi_{\tau^{-n}(\xi)}^{i-1}\circ \phi_{\tau^{-1}(\xi)}^{-n}
(x)\|.
\end{align*}

\begin{claim} \label{cl.dos}
 $s_i(\xi') \leq C_\Phi
(\beta\nu^{\alpha})^{n-i} d_{\Sigma_k}(\xi,\xi')^{\alpha}$,
for every $i \in \{ 1, \dots, n \}$.
\end{claim}

We postpone the proof of this claim.
Taking $n\to \infty$ in \eqref{eq:close-id0}
we get
$$
\|\gamma_{\xi}^{u}(\xi')-\gamma_{\xi}^{u}(\xi)\| \leq  C_\gamma d_{\Sigma_k}(\xi',\xi)^{\alpha} \quad \mbox{for all $\xi' \in W^s_{loc}(\xi;\tau)$ },
$$
where $C_\gamma= C_\Phi(1-\beta\nu^\alpha)^{-1}$,
showing
$\gamma_{\xi}^{u}$ is $\alpha$-H\"older. Indeed for $\xi'$, $\xi'' \in W^u_{loc}(\xi;\tau)$ observe that $\gamma^u_{\xi}(\xi')=\gamma^u_{\xi'}(\xi')$ and $\gamma^u_{\xi}(\xi'')=\gamma^u_{\xi'}(\xi'')$. Thus
$$
  \|\gamma^u_{\xi}(\xi')-\gamma^u_{\xi}(\xi'')\|=
  \|\gamma^u_{\xi'}(\xi')-\gamma^u_{\xi'}(\xi'')\| \leq
  C_\gamma d_{\Sigma_k}(\xi',\xi'')^\alpha.
$$
Note that the H\"older constant obtained is uniform on~$\xi$ and $x=g_\Phi(\xi)$.

\begin{proof}[Proof of Claim~\ref{cl.dos}]
For every $i \in \{ 1, \dots, n\}$, we have that
\[
\begin{split}
s_i(\xi') & =
\|\phi^{n-i}_{\tau^{-n+i}(\xi')}\circ \phi_{\tau^{-n-1+i}(\xi')}
\circ \phi_{\tau^{-n}(\xi)}^{i-1}\circ \phi_{\tau^{-1}(\xi)}^{-n}(x)- \\
& \hspace*{4cm} \phi^{n-i}_{\tau^{-n+i}(\xi')}\circ \phi_{\tau^{-n-1+i}(\xi)}
\circ \phi_{\tau^{-n}(\xi)}^{i-1}\circ \phi_{\tau^{-1}(\xi)}^{-n}(x)\| \\
& \leq   \beta^{n-i}
\| \phi_{\tau^{-n-1+i}(\xi')}\circ \phi_{\tau^{-n}(\xi)}^{i-1}\circ
\phi_{\tau^{-1}(\xi)}^{-n}(x)-  \\
&  \hspace*{4cm}
\phi_{\tau^{-n-1+i}(\xi)} \circ \phi_{\tau^{-n}(\xi)}^{i-1}\circ
\phi_{\tau^{-1}(\xi)}^{-n}(x)\|\\
& \leq   \beta^{n-i}\, C_\Phi \,
d_{\Sigma_k}(\tau^{-n-1+i}(\xi'), \tau^{-n-1+i}(\xi))^\alpha\\
& \leq    (\beta \, \nu^{\alpha})^{n-i}\, C_\Phi \,
d_{\Sigma_k}(\xi', \xi)^{\alpha},
\end{split}
\]
where the first inequality follows from the $\beta$-contraction of $\phi_\xi$, and
the second one from the $\alpha$-H\"older continuity of $\phi_\xi$.
\end{proof}
This completes the proof of Item~\eqref{item-kkk1}.

To prove Item~\eqref{item-kkk2} observe that
\begin{align}
   \phi_{\tau^{-1}(\xi')}^{-1}\circ \gamma_{\xi}^u(\xi')
& = \lim_{n\to \infty} \phi_{\tau^{-1}(\xi')}^{-1} \circ
\phi^{n}_{\tau^{-n}(\xi')}\circ \phi^{-n}_{\tau^{-1}(\xi)}\circ g_\Phi(\xi)\nonumber\\
& = \lim_{n\to \infty} \phi_{\tau^{-1}(\xi')}^{-1} \circ
\phi_{\tau^{-1}(\xi')} \circ
\phi^{n-1}_{\tau^{-n}(\xi')}\circ \phi^{-n}_{\tau^{-1}(\xi)}\circ g_\Phi(\xi)\nonumber\\
& = \lim_{n\to \infty}
\phi^{n-1}_{\tau^{-n}(\xi')}\circ \phi^{-n}_{\tau^{-1}(\xi)}\circ g_\Phi(\xi)\nonumber\\
& = \lim_{n\to \infty} \phi^{n-1}_{\tau^{-n}(\xi')}\circ
\phi^{-(n-1)}_{\tau^{-2}(\xi)}\circ \phi^{-1}_{\tau^{-1}(\xi)}\circ g_\Phi(\xi)\nonumber\\
& = \lim_{n\to \infty} \phi^{n-1}_{\tau^{-n}(\xi')}\circ
\phi^{-(n-1)}_{\tau^{-2}(\xi)}\circ g_\Phi \circ \tau^{-1}(\xi)  \label{e.parakkk2}\\
& = \gamma^u_{\tau^{-1}(\xi)} \circ \tau^{-1}(\xi'), \nonumber
\end{align}
where equality~\eqref{e.parakkk2} is consequence of \eqref{e.tipo conj}:
$\phi^{-1}_{\tau^{-1}(\xi)}\circ
g_\Phi(\xi)= g_\Phi \circ \tau^{-1} (\xi)$.

Now we prove Item~\eqref{item-kkk3}:
$W^{uu}_{loc}((\xi,x);\Phi) \subset W^u((\xi,x);\Phi)$.
Consider $(\xi',x')\in W^{uu}_{loc}((\xi,x);\Phi)$.
Then $\xi' \in W^u_{loc}(\xi;\tau)$ and $x'=\gamma^u_\xi(\xi')$. Thus,
\begin{equation}
\label{eq:ee}
\begin{split}
d(\Phi^{-n}&\big(\xi',x'),\Phi^{-n}(\xi,x)\big) = \\
&=d_{\Sigma_k}(\tau^{-n}(\xi'),\tau^{-n}(\xi)) + \|\phi^{-n}_{\tau^{-1}(\xi')}(x')
- \phi^{-n}_{\tau^{-1}(\xi)}(x) \| \\
&\leq
\nu^n d_{\Sigma_k}(\xi',\xi)+ \|\phi^{-n}_{\tau^{-1}(\xi')}
\circ\gamma^u_{\xi}(\xi')
- \phi^{-n}_{\tau^{-1}(\xi)}(x) \|.
\end{split}
\end{equation}
Since $(\xi,x) \in \Gamma_\Phi$ we have $x=g_\Phi(\xi)$, and using that
$
\gamma^u_{\tau^{-n}(\xi)}\circ \tau^{-n}(\xi) = g_\Phi\circ \tau^{-n}(\xi)=\phi^{-n}_{\tau^{-1}(\xi)}(x) \ \ \text{and} \ \ \phi^{-n}_{\tau^{-1}(\xi')}\circ \gamma^u_{\xi}(\xi') =\gamma^u_{\tau^{-n}(\xi)}\circ \tau^{-n}(\xi')
$
we get 
\[
\begin{split}
\|\phi^{-n}_{\tau^{-1}(\xi')} &\circ\gamma^u_{\xi}(\xi') -
\phi^{-n}_{\tau^{-1}(\xi)}(x) \| = \\
&= \|\gamma^u_{\tau^{-n}(\xi)} \circ\tau^{-n}(\xi')-
\gamma^u_{\tau^{-n}(\xi )}\circ \tau^{-n}(\xi)\|  \leq \nu^{\alpha
n} d_{\Sigma_k}(\xi',\xi)^\alpha.
\end{split}
\]
This implies~\eqref{eq:ee} goes to zero as $n$ goes to infinity and therefore $(\xi',x')$ belongs to $W^u((\xi,x);\Phi)$.  The proof
of Item~\eqref{item-kkk3} and of the proposition is now complete.
\end{proof}

Recall that we denote by $\mathrm{Per}(\Phi)$ the set of periodic
points of $\Phi$ and that $\mathscr{P}: \Sigma_k \times M \to M$
is the projection on the fiber space.

\begin{prop}
\label{p.refe_cube}
 Consider $\Phi=\tau\ltimes \phi_\xi \in \mathcal{S}$
 with $\beta<1$, $\alpha>0$.
 Then
\begin{enumerate}
\addtolength{\itemsep}{.05in}
\item \label{item-ref1}
$ W^{uu}((\xi,x);\Phi)=W^u \big((\xi,x);\Phi \big)\subset
\Gamma_\Phi $ for all  $(\xi,x)\in \Gamma_\Phi$, and
\item \label{item-ref2}
for all periodic point $(\vartheta,p)$ of $\Phi$ in $\Sigma_k
\times \overline{D}$
$$
K_\Phi \eqdef
\overline{\mathscr{P}\big(\mathrm{Per}(\Phi)\big)\cap D }=
\overline{\mathscr{P}\big(W^u((\vartheta,p);\Phi)
\big)}=\mathscr{P}(\Gamma_\Phi) =g_\Phi(\Sigma_k).
$$
\end{enumerate}
\end{prop}

\begin{proof}
To prove the inclusion $W^{uu}((\xi,x);\Phi) \subset
W^{u}((\xi,x);\Phi)$ for all $(\xi,x)$ in $\Gamma_\Phi$, we take
$(\xi',x') \in W^{uu}((\xi,x);\Phi)$ and show that
\begin{equation}
\label{eq:seguida}
   \lim_{n\to \infty} d(\Phi^{-n}(\xi',x'),\Phi^{-n}(\xi,x))=0.
\end{equation}
Since $(\xi',x') \in W^{uu}((\xi,x);\Phi)$, by definition there
are $m\in \mathbb{N}$
and $(\zeta,z) \in W^{uu}_{loc}(\Phi^{-m}(\xi,x);\Phi)$ such that
$(\xi',x')=\Phi^{m}(\zeta,z).$ Let $(\eta,y)=\Phi^{-m}(\xi,x)$.
Note that $(\eta,y)\in \Gamma_\Phi$, $(\zeta,z)\in
W^{uu}_{loc}((\eta,y);\Phi)$, and
$$
   d(\Phi^{-n}(\xi',x'),\Phi^{-n}(\xi,x))=
   d(\Phi^{-(n-m)}(\zeta,z),\Phi^{-(n-m)}(\eta,y)).
$$
By Item~\eqref{item-kkk3} of Proposition~\ref{p:u-lam}, we have
that $W^{uu}_{loc}((\eta,y);\Phi) \subset W^{u}((\eta,y);\Phi)$
and thus it follows~\eqref{eq:seguida}.

Let us prove the converse inclusion: if $(\xi,x)\in \Gamma_\Phi$
then $W^u((\xi,x);\Phi) \subset W^{uu}((\xi,x);\Phi)$. Take
$(\xi',x')\in W^u((\xi,x);\Phi )$ and we have to show that there
is $m\in \mathbb{N}$ such that $\Phi^{-m}(\xi',x')\in
W^{uu}_{loc}(\Phi^{-m}(\xi,x);\Phi)$. By definition of the
unstable set,
\begin{equation}
\begin{split}
\label{e.A} \lim_{n \to \infty} d_{\Sigma_k}
(\tau^{-n}(\xi'),\tau^{-n}(\xi))=0 \quad \text{and} \\
 \lim_{n \to \infty} \|\phi^{-n}_{\tau^{-1}(\xi')}(x')-
\phi^{-n}_{\tau^{-1}(\xi)}(x)\|=0.
\end{split}
\end{equation}
Since $(\xi,x)\in \Gamma_\Phi$ then
$\phi^{-n}_{\tau^{-1}(\xi)}(x)\in D$ for all $n\geq 0$. Thus there
exists $m\in \mathbb{N}$ such that
\begin{equation*}
\tau^{-m}(\xi')\in W^u_{loc} (\tau^{-m}(\xi);\tau ) \quad
\text{and} \quad \phi^{-(n+m)}_{\tau^{-1}(\xi')}(x')\in D, \quad
\mbox{for every $n\geq m$}.
\end{equation*}
Write $(\eta,y)=\Phi^{-m}(\xi,x)\in \Gamma_\Phi$ and
$(\eta',y')=\Phi^{-m}(\xi',x')$. Hence $y=g_\Phi(\eta)$ and
$\phi^{n-i}_{\tau^{-n}(\eta)}\circ\phi^{-n}_{\tau^{-1}(\eta)}(y)$,
$\phi^{n-i}_{\tau^{-n}(\eta')}\circ\phi^{-n}_{\tau^{-1}(\eta')}(y')$
 belong to $D$ for all $0<i \leq n$. Recalling the definition of $\gamma^u_\eta$
and $\gamma^{u,n}_\eta$ we have that
\begin{align*}
\| y'-\gamma^u_{\eta}(\eta')\|&=\lim_{n\to\infty} \|\phi^n_{\tau^{-n}(\eta')}\circ \phi_{\tau^{-1}(\eta')}^{-n}(y')- \phi^{n}_{\tau^{-n}(\eta')}\circ\phi^{-n}_{\tau^{-1}(\eta)}(y) \| \\
&\leq \lim_{n\to\infty} \,\beta^n \,
\|\phi_{\tau^{-1}(\eta')}^{-n}(y')-
\phi^{-n}_{\tau^{-1}(\eta)}(y)\|.
\end{align*}
From~\eqref{e.A} and since $\beta < 1$, we get this limit is zero
and so $y'=\gamma^{u}_{\eta}(\eta')$. That is
$\Phi^{-m}(\xi',x')\in W^{uu}_{loc}(\Phi^{-m}(\xi,x);\Phi)$ and
therefore $(\xi',x')\in W^{uu}((\xi,x);\Phi)$, concluding our
assertion.

Note, we  proved that $ W^{uu}((\xi,x);\Phi)=W^u(\xi,x);\Phi) $
for all $(\xi,x)\in \Gamma_\Phi$. Then by
Proposition~\ref{p:u-lam}, $W^u(\xi,x);\Phi)\subset \Gamma_\Phi$
ending  the first item of the proposition.

To prove Item~\eqref{item-ref2}, consider a periodic point
$\vartheta \in \Sigma_k$ of $\tau$ and note that $W^u(\vartheta;
\tau)$ and $\mathrm{Per}(\tau)$ are both dense in $\Sigma_k$.
Using the conjugation in Proposition~\ref{l.maio2} and
item~\eqref{item-ref1}, we get
\begin{equation}
\label{eq:pp} \overline{\mathrm{Per}(\Phi|_{\Gamma_\Phi})} =
\Gamma_\Phi = \overline{W^u
\big((\vartheta,g_\Phi(\vartheta));\Phi \big)}=\overline{W^{uu}
\big((\vartheta,g_\Phi(\vartheta));\Phi \big)}.
\end{equation}
If $(\vartheta,p)\in \Sigma_k \times \overline{D}$ is a periodic
point of $\Phi$, from the assumption
$\phi_\xi(\overline{D})\subset D$, we have that
$\Phi^n(\vartheta,p) \in \Sigma_k \times D$ for all $n \in
\mathbb{Z}$. Moreover, since $g_\Phi$ is the unique invariant
graph of $\Phi$ restricted to $\Sigma_k \times \overline{D}$, by
Item~\eqref{ite2-thm1} in Theorem~\ref{t.grafo} we get
$p=g_\Phi(\vartheta)$. From this,
\begin{equation}
\label{eq:ppp}
 \mathrm{Per}(\Phi|_{\Gamma_\Phi})=
  \mathrm{Per}(\Phi|_{\Sigma_k\times \overline{D}})
  =\mathrm{Per}(\Phi)\cap (\Sigma_k\times D).
\end{equation}
Recalling that $K_\Phi$ is the closure of projection by
$\mathscr{P}$  of the periodic points of $\Phi$ in $D$ and since
the projection $\mathscr{P}$ is a closed map and $\Sigma_k$ is  a
compact set,
Equations~\eqref{eq:pp} and \eqref{eq:ppp} imply that
\begin{equation*}
\mathscr{P}(\Gamma_\Phi)= \overline{\mathscr{P}
\big(W^u((\vartheta,p);\Phi) \big)} = \overline{\mathscr{P}
\big(\mathrm{Per}(\Phi|_{\Gamma_\Phi}) \big)}
=\overline{\mathscr{P} \big(\mathrm{Per}(\Phi)\big)\cap D }
\eqdef  K_\Phi . 
\end{equation*}
Finally, by definition $\mathscr{P}(\Gamma_\Phi)=g_\Phi(\Sigma_k)$
and from the above equation $K_\Phi =g_\Phi(\Sigma_k)$, completing
the proof.
\end{proof}

Observe that Proposition~\ref{p:u-lam} and~\ref{p.refe_cube}
provide the second and third items in Theorem~\ref{t.B}. Now we
will show that $K_\Phi$ depends continuously with respect to
$\Phi$.
Consider the set $\mathcal{K}(\overline{D})$ whose elements are
the compact subsets of $\overline{D}$ endowed with the Hausdorff
metric
$$
d_H(A,B)\eqdef \sup \{ d(a, B),d(b, A):  a \in A, b \in B \},
\quad \text{$A$, $B \in \mathcal{K}(\overline D)$},
$$
obtaining a complete metric space.
%

Define the map $
 \mathscr{L}:\mathcal{S} \to \mathcal{K}(\overline D)$
given by $ \mathscr{L}(\Phi)\eqdef K_\Phi. $

\begin{prop}
\label{p.cont} Consider $\Phi=\tau \ltimes \phi_\xi \in
\mathcal{S}$ with $\beta <1$, $\alpha\geq 0$. Then, for each
$\varepsilon>0$ there is $\delta>0$ such that for every
$\Psi=\tau\ltimes \psi_\xi \in \mathcal{S}$ with
$$
d_{C^0}(\phi_\xi,\psi_\xi)<\delta \quad \text{it holds  that}
\quad d_H (K_\Phi, K_\Psi ) = d_H (\mathscr{L} (\Phi),
\mathscr{L}(\Psi))
 <\varepsilon.
$$
\end{prop}

The above result follows from Proposition~\ref{t.grafo} and the
inequality
$$
d_H (K_\Phi, K_\Psi) = d_H \big( g_\Phi(\Sigma_k),
g_\Psi(\Sigma_k) \big) \leq d_{C^0} ( g_\Phi, g_\Psi).
$$
This completes the proof of Theorem~\ref{t.B}.

\section*{Acknowledgements}
We are grateful to Lorenzo J. D\'iaz  for his assistance and review of preliminar versions of this paper and
Enrique R. Pujals for helpful conversations and comments.
The first two authors would also like to express a deep gratitude to their respective scientific advisors, J. \'Angel Rodr\'iguez and Lorenzo J. D\'iaz for the support.

During the preparation of this article the authors were supported by the following fellowships:  P. G. Barrientos by FPU-grant AP2007-031035,  MTM2011-22956 project (Spain) and CNPq post-doctoral fellowship (Brazil);
Y. Ki by CNPq-Brazil doctoral fellowship and CNPq P\'os doutorado J\'unior 503802/2012-3; A. Raibekas by CNPq-Brazil doctoral and CAPES-PNPD-Brazil post-doctoral fellowships.

We thank Instituto Nacional de Matem\'atica Pura e Aplicada (IMPA)-Brazil for its hospitality.

\bibliographystyle{alpha}
\bibliography{bkr}

\end{document}